\documentclass[a4paper,11pt]{amsart}

\usepackage{amsmath,amssymb,amsthm}
\usepackage{amsfonts}
\usepackage{pgf,tikz}
\usepackage{paralist}
\usepackage{graphicx}
\usepackage{float}
\usepackage{caption}
\usepackage{subcaption}
\usepackage{amsrefs}
\usepackage{mathrsfs}

\usetikzlibrary{arrows}
\usepackage{color}

\definecolor{light-gray}{gray}{0.85}
\usepackage[colorlinks=true,linkcolor=blue,urlcolor=blue]{hyperref}
\usepackage{listings}
\usepackage[colorlinks=true,linkcolor=blue,urlcolor=blue]{hyperref}
\usepackage{aliascnt}

\theoremstyle{plain}
\newtheorem{theorem}{Theorem}[section]

\newaliascnt{cor}{theorem}
\newtheorem{cor}[cor]{Corollary}
\aliascntresetthe{cor}

\newaliascnt{lem}{theorem}
\newtheorem{lem}[lem]{Lemma}
\aliascntresetthe{lem}

\newaliascnt{conjecture}{theorem}

\aliascntresetthe{conjecture}

\newaliascnt{prop}{theorem}
\newtheorem{prop}[prop]{Proposition}
\aliascntresetthe{prop}

\newaliascnt{dfn}{theorem}
\newtheorem{dfn}[dfn]{Definition}
\aliascntresetthe{dfn}

\newaliascnt{nota}{theorem}
\newtheorem{nota}[nota]{Notation}
\aliascntresetthe{nota}

\newaliascnt{exmp}{theorem}
\newtheorem{exmp}[exmp]{Example}
\aliascntresetthe{exmp}

\newaliascnt{rem}{theorem}
\newtheorem{rem}[rem]{Remark}
\aliascntresetthe{rem}

\newaliascnt{observe}{theorem}
\newtheorem{observe}[observe]{Observation}
\aliascntresetthe{observe}

\newaliascnt{algorithm}{theorem}
\newtheorem{algorithm}[algorithm]{Algorithm}
\aliascntresetthe{algorithm}

\newtheorem{headtheorem}{Theorem}

\newcommand{\hmatTagSize}[2]{\mathbf{X}^{#1}_{#2}}

\newcommand{\rees}[1]{\mathscr{R}(#1)}
\newcommand{\reesSymbol}{\mathscr{R}}
\newcommand{\reesInitial}[1]{\mathscr{R}^{\operatorname{in}}(#1)}
\newcommand{\kerRees}{\ker_{\reesSymbol}}
\newcommand{\kerReesIn}{\ker_{\reesSymbol^{\operatorname{in}}}}
\newcommand{\initial}[1]{\operatorname{in}(#1)}
\newcommand{\initiaL}{\operatorname{in}}
\newcommand{\chain}[2]{{#1}_1,\ldots,#1_{#2}}
\newcommand{\f}{\mathcal{F}}

\newcommand{\diagonalSquence}[2]{(#1,#2)}
\newcommand{\diagSeqA}{(\sigma,a)}
\newcommand{\diagSeqC}{(\sigma,c)}
\newcommand{\diagSeqB}{(\gamma,b)}
\newcommand{\diagSeqD}{(\gamma,d)}
\newcommand{\tabelA}[2]{\mathcal{A}=\{\diagonalSquence{#1_i}{#2^{(i)}}\}}
\newcommand{\tabelB}[2]{\mathcal{B}=\{\diagonalSquence{#1_i}{#2^{(i)}}\}}
\newcommand{\interval}[1]{\mathscr{L}(#1)}

\newcommand{\detIdeal}[2]{I_{#1,#2}}

\usepackage{xcolor}
\usepackage[normalem]{ulem} % use normalem to protect \emph
\newcommand\hl{\bgroup\markoverwith
	{\textcolor{yellow}{\rule[-.5ex]{2pt}{2.5ex}}}\ULon}

\title{Linear resolutions and Gr\"{o}bner basis of Hankel determinantal ideals}

\author{Sepehr Jafari}
\address {Dipartimento di Matematica, Università degli studi di Genova, Via Dodecaneso, 35, 16146 Genova GE, Italy}
\email{sepehr@dima.unige.it}



\begin{document}

\maketitle

\begin{abstract}
	In this paper, we study the family of determinantal ideals of "close" cuts of Hankel matrices, say $ \f $. We show that the multi-Rees algebra of ideals in $ \f $ is defined by a quadratic Gr\"{o}bner basis, it is Koszul, normal Cohen-Macaulay domain  and it has a nice Sagbi basis. As a consequence of Koszulness, we prove that every product $ I_1\ldots I_l $ of ideals of $ \f $ has linear resolution. Moreover, we show that natural generators of every  product  $ I_1\ldots I_l $ form a Gr\"{o}bner basis. 
\end{abstract}

\section*{Introduction}

The study of determinantal  ideals is a classic topic in commutative algebra. The main properties of determinantal ideals are described in the book by Bruns and Vetter (see   \cite{bruns-vetter}). The cases of generic symmetric and generic skew-symmetric matrices, are also well-understood thanks to  J\'ozefiak and Pragacz (see \cite{generic_symmetric} and \cite{generic_skew_symmetric}). One standard method to study determinantal ideals is to understand their initial   ideals via Gr\"{o}bner basis. Often, Gr\"{o}bner basis of determinantal ideals nicely encodes the information of the minors (i.e the product of the entries on the diagonal or anti-diagonal). For the  generic matrices Sturmfels and separately Herzog and  Trung (see \cite{G_basis_Sturmfels} and \cite{Herzog_Trung_G_basis_multiplicity}) described the Gr\"{o}bner basis. As for symmetric matrices and in particular Hankel matrices, the Gr\"{o}bner basis is given by Conca  (see \cite{Conca_G_basis_symmetric_matrix} and \cite{CONCA_STR_LAW}).\\
The question of the study of Castelnuovo-Mumford regularity of products and powers of ideals has attracted many researchers in the past decades. One important result in this direction was given in 1999 by  Cutkosky, Herzog and Trung (see \cite{Cutkosky}) and independently in 2000 by  Kodiyalam (see \cite{Kodiyalam}). The result, in essence, is that the regularity of large enough powers of  ideal $ I $ in  polynomial ring $ S $ is given by a linear function. We say ideal $ I $ generated in degree $ d $ has linear resolution if and only if regularity of $ I $ coincides with $ d $. It is interesting to find family of ideals with (asymptotic) linear powers  (i.e The family $ \f $ of ideals of polynomial ring $ S $  has (asymptotic) linear powers if $ I^q $ has linear resolution for  every  $ I\in\f$ and (large enough) every positive number $ q $.). The literature on finding families of ideals with (asymptotic) linear powers is vast. In 2015, Bruns, Conca and Varbaro proved that the determinantal ideals of maximal minors of a generic matrix has linear powers    (see \cite{Bruns_Conca_Varbaro}). In 2018, Raicu classified the determinantal  ideals of a generic matrix with asymptotic linear powers (see \cite{Raicu}).\\ 
The asymptotic behavior of Betti numbers of products of ideals is studied by Bagheri, Chardin and Ha in 2013 (see \cite{Bagheri_Chardin_Huy}). In 2017, Bruns and Conca described the asymptotic behavior of regularity of products of ideals with a different method (see \cite{Bruns_Conca_A_remark_regularity}). 
Of particular interest is to find families of ideals with linear products (i.e The family $ \f $ of  ideals of  polynomial ring $ S $  has linear products if every product $ I_1\ldots I_l $  of the ideals in $ \f $ has linear resolution). In 2003,  Conca and    Herzog introduced few families of ideals with linear products. In particular, the authors showed that the family of determinantal ideals of Hankel matrices with entries indeterminates $ \chain{x}{n} $ of ring $ S $  has linear products (see \cite{Conca_Herzog_reg_products_ideals}). In 2011,  Nam improved the former result to the family of extended Hankel matrices (see  \cite{NAM_PAPER}). For generic matrices, in 2017, Bruns and Conca proved that  "north-east" determinantal ideals of a generic matrix (or Borel fixed ideals as the authors call them.) forms a family of ideals with linear products (see   \cite{Bruns_Conca}). Moreover, in the above works (\cite{Conca_Herzog_reg_products_ideals}, \cite{NAM_PAPER} and \cite{Bruns_Conca} ), the authors prove as a consequence that  every product of ideals, in their cases of study, has a Gr\"{o}bner basis given by its natural generators with respect to some term order.\\
Inspired by Bruns and Conca \cite{Bruns_Conca}, we introduce a new family of determinantal ideals of Hankel matrices with linear products and nice  Gr\"{o}bner basis to improve the result of Nam in \cite{NAM_PAPER}. We take advantage of Sagbi deformations theory. The theory of Sagbi basis (i.e one can see it as analogues of Gr\"{o}bner basis for algebras) was introduced by Robbiano and Sweedler (see \cite{Robbiano_Sweedler}) and independently by  Kapur and Madlener (see \cite{Kapur_Madlener}). In 1996, Conca, Herzog and Valla introduced Sagbi deformations and applied it  to study Rees algebras of certain rational normal scrolls (see \cite{Conca_Herzog_Valla}).

Let $ S=K[\chain{x}{n}] $ be a standard graded polynomial ring with coefficients from arbitrary field $ K $.  Let  $ \hmatTagSize{(1,n)}{t} $, where $ 1\leq t \leq \lfloor\frac{n+1}{2}\rfloor $, denote the Hankel matrix with $ t $ rows and entries $ x_1,\ldots,x_n $:
\begin{center}
	$ \hmatTagSize{(1,n)}{t}=\begin{pmatrix}
	x_{1}       & x_{2} & x_{3} & \dots & x_{n-t+1} \\
	x_{2}       & x_{3} & \dots & \dots & \dots \\
	x_{3}       & \dots & \dots & \dots & \dots \\
	\vdots       & \vdots & \vdots & \vdots & \vdots \\
	x_{t}       & \dots &\dots & \dots & x_{n}
	\end{pmatrix} $
\end{center}

Let $ \hmatTagSize{(1,n)}{} $ be the family of all Hankel matrices with entries $ \chain{x}{n} $. Denote  by $ \f^{(1,n)} $ the family of all determinantal ideals of matrices in $ \hmatTagSize{(1,n)}{} $. Let $ \detIdeal{(1,n)}{t}  $ denote the determinantal ideal of maximal minors of $ \hmatTagSize{(1,n)}{t} $. It is known that $ \detIdeal{(1,n)}{2} $ defines the rational normal curve as well as $ \detIdeal{(1,n)}{t} $ defines the $ t- $th secant variety of the rational normal curve, moreover, the minimal free resolution of $ \detIdeal{(1,n)}{t} $ is known to be the Eagon-Northcott resolution.\\
Let $ \hmatTagSize{(i,j)}{} $ and $ \f^{(i,j)} $ with $ i+1<j $ be defined as above. By experiments, we expect $ \cup_{i+1<j}\f^{(i,j)} $ to have linear products. It is not difficult to see that the Sagbi deformation fails for $ \cup_{i+1<j}\f^{(i,j)} $ in general (see \autoref{importance}). However,  we can still employ Sagbi deformations for $ \f=\f^{(1,n)}\cup\f^{(1,n-1)}\cup\f^{(2,n)}\cup\f^{(2,n-1)} $. We refer to the family $ \hmatTagSize{(1,n)}{}\cup\hmatTagSize{(1,n-1)}{}\cup\hmatTagSize{(2,n)}{}\cup \hmatTagSize{(2,n-1)}{} $ by "close" cuts of Hankel matrices or shortly close cuts.\\
\noindent In  Section \ref{Standard_Forms}, we study the products of the monomials laying on the diagonal of maximal minors of Hankel matrices. We encode these monomials in a tabel. Then, we transform these tabels into so called standard forms. This machinery is the foundation of our treatments in this paper.\\
In  Section \ref{Multi-Rees_Algebra}, we apply the tools of Gr\"{o}bner basis and  standard forms to study the multi-Rees algebra $ \reesSymbol^{\operatorname{in}}=\rees{\initial{\detIdeal{\sigma}{a}}:\detIdeal{\sigma}{a}\in\f} $. In particular, we prove that $ \reesSymbol^{\operatorname{in}} $ is defined by a quadratic Gr\"{o}bner basis and it is Koszul (see \autoref{main_theorem_initial}). Then, we "lift" this property to the multi-Rees algebra $ \rees{\detIdeal{\sigma}{a}:\detIdeal{\sigma}{a}\in\f} $ by applying  Sagbi deformation. 
Our main result is the following theorem:
\begin{headtheorem}[\autoref{main_theorem}]
	The family $ \f $ has the following features:
	\begin{enumerate}[(1)]
		\item Every product $ \prod_{(\sigma,a)}\detIdeal{\sigma}{a} $ of ideals in $ \f $ has linear resolution.
		\item Computing the initial ideals commutes over products  $ \initial{\prod_{(\sigma,a)}\detIdeal{\sigma}{a}}=\prod_{(\sigma,a)}\initial{\detIdeal{\sigma}{a}} $, in particular the natural generators form a Gr\"{o}bner basis. 
		\item The multi-Rees algebra $ \rees{\detIdeal{\sigma}{a}:\detIdeal{\sigma}{a}\in\f} $ is defined by a quadratic Gr\"{o}bner basis, it is Koszul,  normal, Cohen-Macaulay domain. Moreover, the natural algebra generators form a Sagbi basis.
	\end{enumerate}  
\end{headtheorem}

\noindent In this work, we strongly took advantage of several computations made by computer algebra systems Macaulay2 \cite{M2} and Cocoa5 \cite{CoCoA}.  

\section*{Acknowledgment}
The author would like to show his deep gratitude to his PhD thesis advisor Aldo Conca for his mentorship, suggesting this problem and several helpful discussions.

\section{Labeled chains and Standard Forms}\label{Standard_Forms}

Let $ S=K[x_1,\ldots,x_n] $ be the polynomial ring  equipped with standard grading and degree lexicographical term order with respect to $ x_1>\ldots>x_n $. We denote the term order of $ S $ by $ \leq_{\tau} $. By $ \hmatTagSize{(1,n)}{t} $, where $ 1\leq t \leq \lfloor\frac{n+1}{2}\rfloor $, denote the Hankel matrix with $ t $ rows and entries $ x_1,\ldots,x_n $:
\begin{center}
	$ \hmatTagSize{(1,n)}{t}=\begin{pmatrix}
		x_{1}       & x_{2} & x_{3} & \dots & x_{n-t+1} \\
		x_{2}       & x_{3} & \dots & \dots & \dots \\
		x_{3}       & \dots & \dots & \dots & \dots \\
		\vdots       & \vdots & \vdots & \vdots & \vdots \\
		x_{t}       & \dots &\dots & \dots & x_{n}
	\end{pmatrix} $
\end{center}
For $ s\leq t $, the standard notation for minors of $ \hmatTagSize{(1,n)}{t} $ is $ [\chain{a}{s}|\chain{b}{s}] $ where $ \chain{a}{s} $ and  $ \chain{b}{s} $ are strictly increasing sequences of row indices and column indices respectively. A minor of the form $ [1,\ldots,t|\chain{b}{t}] $ is called a maximal minor of $ \hmatTagSize{(1,n)}{t} $. Unless otherwise is stated, $ \initial{f} $ will always denote the leading term of the polynomial $ f $. Similarly, $ \initial{I} $ will always denote the leading term ideal of the ideal  $ I $. It is clear that the leading term of a minor is  the product of the the entries laying on the main diagonal. The term orders satisfying this criteria are known as diagonal term orders in the literature. The arguments of this paper holds identically for any given diagonal term order. Therefore, our results are true for any given diagonal term order. Let $ i $ and $ j $ be distinct natural numbers. We define the partial order $ \leq_1 $ by  
\begin{center}
	$ i\leq_1 j $ if and only if $ i+1\leq j $.
\end{center}
When $ i+1< j $, we denote the above partial order by $ <_1 $. A chain  is a  sequence in $ \{1,\ldots ,n\} $ like $ a=\chain{a}{r} $ such that $ a_1<_1 a_{2} <_1  \ldots <_1  a_r $. Similarly, we say a monomial $ x_a=x_{a_1}\ldots x_{a_r} $ is a chain if its indices form a chain. A given monomial $ x_{a_1}\ldots x_{a_r} $ is a chain if and only if  it is the leading term of a minor of some $ \hmatTagSize{(1,n)}{t} $ with $ r \leq t $. This corresponding minor is unique if  and only if $ r=t $. We will denote the family of all Hankel matrices with entries $ x_1,\ldots,x_n $ by $ \hmatTagSize{(1,n)}{} $. 
 Let $ \detIdeal{(1,n)}{t} $ denote the ideal generated by the maximal minors of $ \hmatTagSize{(1,n)}{t} $. We denote the family of all determinantal ideals of the matrices of the family $ \hmatTagSize{(1,n)}{} $ by $ \f^{(1,n)} $. It is clear that one can repeat the same constructions for the sequences of indeterminates $ x_1,\ldots,x_{n-1} $, $ x_2,\ldots,x_{n} $  and $ x_2,\ldots,x_{n-1} $ and construct the families of matrices $ \hmatTagSize{(1,n-1)}{} $, $ \hmatTagSize{(2,n)}{} $ and $ \hmatTagSize{(2,n-1)}{} $ and of course $ \f^{(1,n-1)} $, $ \f^{(2,n)} $ and $ \f^{(2,n-1)} $. We will refer to the family 
 $$ \hmatTagSize{}{}=\hmatTagSize{(1,n)}{}\cup\hmatTagSize{(1,n-1)}{}\cup\hmatTagSize{(2,n)}{}\cup\hmatTagSize{(2,n-1)}{} $$ 
 by the family of close Hankel matrices or shortly close cuts. In this work, we investigate the following family:
 $$ \f=\f^{(1,n)}\cup\f^{(1,n-1)}\cup\f^{(2,n)}\cup\f^{(2,n-1)} $$
 We will refer to the tuples $ (1,n),(1,n-1),(2,n) $ and $ (2,n-1) $ by labels. To keep a simple notation, we choose a general notation for them like: 
 $$ \sigma=(\sigma_1,\sigma_2)\in\{(1,n),(1,n-1),(2,n),  (2,n-1)\} $$
 Given a chain $ a=\chain{a}{r} $, there is at least one label, say $ \sigma $, such that the matrix $ \hmatTagSize{\sigma}{r} $ contains a maximal minor whose initial term is $ x_a $. For many cases,  this label is not unique.  A labeled chain $ (\sigma,a) $ is a chain together with a fixed label. For a given labeled chain, we keep the  label fixed  unless otherwise is clearly stated. \\
In our treatment, we need to put an order on the labeled chains. To this end, we order the set of labels lexicographically  like the following: 
$$ (1,n-1)>(1,n)>(2,n-1)>(2,n) $$
Let $ \diagSeqA $ and $ \diagSeqB $ be labeled chains. We say $ \diagSeqA\geq_c \diagSeqB $ if and only if $ \sigma>\gamma $ or $ \sigma=\gamma $ and $ a>_{\tau}b $. We denote a pair of labeled chains by $ \diagSeqA\geq_c \diagSeqB $ when we need to  emphasis on their order. 

\begin{dfn}
	Let $ (\sigma_1,a^{(1)}),(\sigma_2,a^{(2)}),\ldots,(\sigma_k,a^{(k)}) $ be a set of labeled chains where $ a^{(i)}=a^{(i)}_{1},\ldots,a^{(i)}_{r_i} $ for all $ i\colon1,\ldots,l $.	A tabel $ \tabelA{\sigma}{a} $ is a set of labeled chains such that  $ (\sigma_i,a^{(i)})\geq_c(\sigma_{i+1},a^{(i+1)}) $ for all $ i\colon1,\ldots,l-1 $. We will refer to $ (\sigma_i,a^{(i)}) $ by the $ i-$th row of $ \mathcal{A} $.
\end{dfn}

\begin{nota}
	Throughout this paper, we reserve $ r_i $ to denote the length of the $ i-$th row of some tabel $ \mathcal{A} $. The symbol $ \sigma_i $ is reserved for the notation of the  label of the $ i-$th row of $ \mathcal{A} $. We reserve the letter $ l $ to denote the number of the rows of $ \mathcal{A} $. In particular, when  we mention chains $ a $ and $ b $, we denote the length of $ a $ by $ r $ and the one of $ b $ by $ s $. 
\end{nota}

\begin{dfn}
	Let $ \tabelA{\sigma}{a} $ be a tabel. The shape of $ \mathcal{A} $ is the sequence $ (\sigma_1,r_1),(\sigma_2,r_2),\ldots,(\sigma_l,r_l) $. 
\end{dfn}

Recall that a tabel $ \tabelA{\sigma}{a} $ contains the information on the monomial $ \prod_{\mathcal{A}}x_{a_i} $ together with the supporting labels. In the treatment of this section, we need to transform $ \mathcal{A} $ to an other tabel, say $ \mathcal{B} $, in such a way that both tabels encode the same monomial and have the same shape. The following function $ \Delta $ controls the shape of $ \mathcal{A} $ during this  transformation. The functions $ \Omega_{\text{d}} $, $ \Omega_{\text{c}} $ and $ \Omega_{\text{ad}} $ control the iteration through the transformations of pairs of rows. 

\begin{dfn}\label{controllers}
	Let  $ (\sigma,a) \geq_c (\gamma,b) $  be a pair of labeled chains. We define 
	\begin{align*}
		&\Delta\{\diagSeqA\geq_c \diagSeqB \} =\begin{cases}
		1 &\text{if } \quad r\leq s, a_r=n \text{ and }\gamma_2=n-1,\\
		&\text{or} \quad r> s, b_s=n\text{ and }\sigma_2=n-1,\\
		0  &\text{otherwise}.
		\end{cases}\\
		&\Omega_{d}\{\diagSeqA\geq_c \diagSeqB\} =\begin{cases}
		r-\Delta\{\diagSeqA\geq_c \diagSeqB\} &\text{ if } \quad r< s,\\
		s-1  &\text{ if } \quad r\geq s.
		\end{cases}\\
		&\Omega_{c}\{\diagSeqA\geq_c \diagSeqB\} =\begin{cases}
		r-\Delta\{\diagSeqA\geq_c \diagSeqB\}  &\text{ if } \quad r\leq s,\\
		s & \text{ if } \quad r>s.
		\end{cases}\\
		&\Omega_{ad}\{\diagSeqA\geq_c \diagSeqB\} =\begin{cases}
		r-1 &\text{ if } \quad r\leq s,\\
		s-\Delta\{\diagSeqA\geq_c \diagSeqB\}  &\text{ if } \quad r> s.
		\end{cases}
	\end{align*}

\noindent In case there is no confusion, we use $ \Delta,\Omega_{d}, \Omega_{c} $ and $ \Omega_{ad} $ instead of  $ \Delta\{(\sigma,a),(\gamma,b)\} $,  $ \Omega_{d}\{(\sigma,a),(\gamma,b)\}$, $ \Omega_{c}\{(\sigma,a),(\gamma,b)\} $ and $ \Omega_{ad}\{(\sigma,a),(\gamma,b)\} $ respectively. We present these functions in the following tabel: 
\begin{center}
	\begin{tabular}{ c|c|c|c }\label{tabel} 
		& $ \Omega_d $ & $ \Omega_c $ & $ \Omega_{ad} $ \\ 
		\hline
		$ r<s $ & $ r-\Delta $ & $ r-\Delta $ & $ r-1 $ \\ 
		\hline
		$ r>s $ & $ s-1 $ & $ s $ & $ s-\Delta $ \\ 
		\hline
		$ r=s $ & $ s-1 $ & $ s-\Delta $ & $ s-1 $ \\ 	
		\end{tabular}	
\end{center}
\end{dfn}

Let $a=\chain{a}{r}$ be a chain. We define $\interval{a}=\bigcup_{2\leq i \leq r}\{a_i-1,a_i\}$.
Let $ \diagSeqA\geq_c \diagSeqB $ be a pair labeled chains. We have the following relations: 
\begin{itemize}
	\item [] \textbf{Diagonal relation}: The pair  $ \diagSeqA\geq_c \diagSeqB $ has diagonal relations if $ a_{h}>b_{k} $, for some $ 1\leq h \leq \Omega_{\text{d}} $ and $ h+1\leq k \leq s $ and $ a_{h}\notin \interval{b} $. If $ \diagSeqA\geq_c \diagSeqB $ does not have any diagonal relations, we say it is diagonal sorted.
	\item [] \textbf{Column-wise relations}: The pair  $ \diagSeqA\geq_c \diagSeqB $ has column-wise relations if it  is diagonal sorted and $ a_{h}>b_{h} $ for some $ 1\leq h \leq  \Omega_{\text{c}} $ and  $ a_{h}\notin \interval{b} $. If $ \diagSeqA\geq_c \diagSeqB $ does not have any column-wise relations, we say it is column-wise sorted.
	\item [] \textbf{Anti Diagonal relations}: The pair  $  \diagSeqA\geq_c \diagSeqB $ has anti-diagonal relations if it is column-wise sorted and $ b_{h}>a_{k} $ for some $ 1\leq h \leq  \Omega_{\text{ad}} $ and $ h+1\leq k \leq r $ and $ b_{h}\notin \interval{a} $. If $ \diagSeqA\geq_c \diagSeqB $ does not have any anti-diagonal relations, we say it is anti-diagonal sorted.
\end{itemize}
With respect to the above relations, we define what we consider as the standard form.
\begin{dfn}\label{std_table}
	The pair $ \diagSeqA\geq_c \diagSeqB $ is a standard form if it satisfies the following conditions: 
	\begin{itemize}
		\item [(1)] $ a_i\leq b_{i+1} $ for $ i\colon 1,\ldots,\Omega_{\text{d}} $ or  $ a_h>b_{h+1} $ for some $ h $ and $ a_h,\ldots,a_{\Omega_{\text{d}}}\in \interval{b} $;
		\item [(2)] $ a_i\leq b_{i} $ for $ i\colon 1,\ldots,\Omega_{\text{c}} $ or  $ a_h>b_{h} $ for some $ h $ and $ a_h,\ldots,a_{\Omega_{\text{c}}}\in \interval{b} $;
		\item [(3)]  $ b_{i}\leq a_{i+1}  $ for $ i\colon 1,\ldots,\Omega_{\text{ad}} $ or  $ b_h>a_{h+1} $ for some $ h $ and $ b_h,\ldots,b_{\Omega_{\text{ad}}}\in \interval{a} $.
	\end{itemize}
	We say a tabel $ \tabelA{\sigma}{a} $ is a standard form if the following conditions hold for $ 1\leq s < t \leq k $:
	\begin{itemize}
		\item [(1)] $ a^{(s)}_{i}\leq a^{(t)}_{i+1} $ for $ i\colon 1,\ldots,\Omega_{\text{d}} $ or  $ a^{(s)}_{h}>a^{(t)}_{h+1} $ for some $ h $ and $ a^{(s)}_{h},\ldots,a^{(s)}_{\Omega_{\text{d}}}\in \interval{a^{(t)}} $;
		
		\item[(2)] $ a^{(s)}_{i}\leq a^{(t)}_{i} $ for $ i\colon 1,\ldots,\Omega_{\text{c}} $ or  $ a^{(s)}_{h}>a^{(t)}_{h} $ for some $ h $ and $ a^{(s)}_{h},\ldots,a^{(s)}_{\Omega_{\text{c}}}\in \interval{a^{(t)}} $;
		
		\item [(3)]  $ a^{(t)}_{i}\leq a^{(s)}_{i+1} $ for $ i\colon 1,\ldots,\Omega_{\text{ad}} $ or  $ a^{(t)}_{h}>a^{(s)}_{h+1} $ for some $ h $ and $ a^{(t)}_{h},\ldots,a^{(t)}_{\Omega_{\text{ad}}}\in \interval{a^{(s)}} $.
	\end{itemize}
\end{dfn}

In the second part of the   \autoref{std_table}, we refer to a tabel satisfying  $ (1) $, $ (2) $ and $ (3) $ by diagonal sorted, column-wise sorted and anti-diagonal sorted respectively.\\

\begin{prop}
	A pair of labeled chains  $ \diagSeqA\geq_c \diagSeqB $  is a standard form if and only if it is diagonal sorted, column-wise sorted and anti-diagonal sorted. 
	\begin{proof}
		If $  \diagSeqA\geq_c \diagSeqB  $ is a standard form, it is clear that it is diagonal sorted, column-wise sorted and anti-diagonal sorted. It remains to prove the other direction. \\
		Let $ d=\Omega_d $. We show \autoref{std_table} $ (1) $ holds.  Suppose $ a_h>b_{h+1} $ for some $ 1\leq h \leq d $. The case $ h=d $ is trivial. Since the pair $ \diagSeqA\geq_c \diagSeqB $ is reduced modulo the relations, we have $ a_h\in \interval{b} $. Thus, there exists unique $ h+1 < t \leq s $  such that $ a_h\in\{b_{t}-1,b_{t}\} $. Since $ a_{h}<_1 a_{h+1} $ we have $ a_h \leq b_{t} <a_{h+1} $. Looking at the bounds of $ t $ and the fact that $ \diagSeqA\geq_c \diagSeqB $ is  is diagonal sorted, we have $ a_{h+1}\in\interval{b} $. By repeating this argument, we obtain $ a_{h+1},\ldots,a_d\in\interval{b} $. Therefore the first part of the definition holds. \\
		Let $ d=\Omega_c $. We show \autoref{std_table} $ (2) $ holds. Suppose $ a_h>b_h $ for some $ 1\leq h \leq d $. The case $ h=d $ is trivial. Since $ \diagSeqA\geq_c \diagSeqB $ is reduced modulo the relations, we have $ a_h\in\interval{b} $. So there exists $ h < t \leq s $ such that $ a_h\in\{b_{t}-1,b_{t}\} $. From $ a_h<_1 a_{h+1} $, we have  $ a_h\leq b_{t} < a_{h+1} $. Since $ \diagSeqA\geq_c \diagSeqB $ is reduced modulo the relations (in particular it is diagonal sorted and column-wise sorted), we have $ a_{h+1}\in\interval{b} $. Repeating this argument gives $ a_{h+1},\ldots,a_d\in\interval{b} $.\\ 
		Let $ d=\Omega_{ad} $. We show \autoref{std_table} $ (3) $ holds. Suppose $ b_h>a_{h+1} $ for some $ 1\leq h \leq d $. The case $ h=d $ is trivial.  Since $ \diagSeqA\geq_c \diagSeqB $ is reduced modulo relations, we have $ b_h\in\interval{a} $ which gives a unique $ h+1 < t \leq r $ such that $ b_h\in\{a_{t}-1,a_{t}\} $. From $ b_h<_1b_{h+1} $, we get $ b_h\leq a_{t} < b_{h+1} $. Since $ \diagSeqA\geq_c \diagSeqB $ is reduced modulo relations (in particular it is anti-diagonal sorted), we have $ b_{h+1}\in\interval{a} $. By repeating this argument, we get $ b_h,\ldots,b_d\in\interval{a} $.
	\end{proof}
\end{prop}

\begin{rem}\label{largest_row_last_entry}
	The following are easy to check:
	\begin{enumerate} [(1)]
		\item  The tabel $ \tabelA{\sigma}{a} $ is a standard form if and only if the pair $ (\sigma_i,a^{(i)})\geq_c(\sigma_j,a^{(j)}) $ is a standard form for all $ 1\leq i < j \leq k $.
		\item If the tabel $ \tabelA{\sigma}{a} $ is a standard form and $ (\sigma_i,a^{(i)}) $ and $ (\sigma_j,a^{(j)}) $ are some rows of $ \mathcal{A} $, then
		\begin{enumerate}[(i)]
			\item If    $ r_i< r_j $, when $ a^{(i)}_{r_i}\leq \sigma_{j,2} $ we have $ a^{(i)}_{r_i}\leq a^{(j)}_{r_j} $, otherwise we have $ a^{(i)}_{r_i}> a^{(j)}_{r_j} $.
			\item If   $ r_i=r_j $ and $ i<j $, when $ \Delta=0 $ we have $ a^{(i)}_{r_i}\leq a^{(j)}_{r_j} $, otherwise $ a^{(i)}_{r_i}> a^{(j)}_{r_j} $.
		\end{enumerate}
	\end{enumerate}
\end{rem}

\begin{lem}\label{diagonal_reduction}
	Any pair of labeled chains  $ \diagSeqA\geq_c \diagSeqB $ transforms to a diagonal sorted   form $ \diagSeqC \geq_c \diagSeqD $ of the same shape. 
	\begin{proof}
		Suppose $ \diagSeqA\geq_c \diagSeqB $ is not diagonal sorted. There exists $ 1\leq h \leq \Omega_d $ and $ h+1\leq k \leq s $ such that $ a_h>b_k $ and $ a_h\notin\interval{b} $. Assume $ (h,k) $ is minimum with respect $ \leq_\tau $. Recall that $ \leq_\tau $ is $ \operatorname{lex} $ order induced by $ x_1>x_2>\ldots>x_n $.We prove the statement by  induction on $ (h,k) $. There exists $ 0\leq v \leq h-1 $ such that for every $ 0\leq i \leq v $, we have $ a_{h-i}>b_{k-i} $ and $ a_{h-v-1}\leq b_{k-v-1} $ if $ v\neq h-1 $. Clearly, for $ v \neq h-1 $, we have $ a_{h-v-1}\leq b_{k-v-1} <_1 b_{h-v} $ and $ b_{k}< a_{h}<_1 a_{h+1} $. Therefore, $ \tilde{a}=\chain{a}{h-v-1},b_{k-v},\ldots,b_{k},a_{h+1},\ldots,a_r $ is a chain. By definition of diagonal relations and $ \Omega_d $, one can see that $ (\sigma , \tilde{a}) $ is a well-defined labeled chain. On the other hand, $ b_{k-v-1}<_1 b_{k-v}<a_{h-v} $ and $ a_h<_1 b_{k+1} $ from  $ a_h\notin\interval{b} $ and the definition of $ k $. Therefore, $ \tilde{b}=b_1,\ldots,b_{k-v-1},a_{h-v},\ldots,a_{h},b_{k+1},\ldots,b_{s} $ is a chain. By definition of diagonal relations and $ \Omega_d $, one can see that $ (\gamma,\tilde{b}) $ is well-defined. Moreover, from the above and the definition of $ \geq_c $, one can see $ (\sigma , \tilde{a})\geq_c(\gamma,\tilde{b}) $.  Let $ (\tilde{h},\tilde{k}) $ be the analogous of $ (h,k) $ for the new pair  $ (\sigma , \tilde{a})\geq_c(\gamma,\tilde{b}) $. It remains to show that $ (\tilde{h},\tilde{k})\geq_\tau(h,k) $, hence the induction yields.  By  construction of $ \tilde{a} $ and $ \tilde{b} $, we have $ \tilde{a}_{h+1},\tilde{a}_{h+2},\ldots,\tilde{a}_{r}=a_{h+1},a_{h+2},\ldots,a_r $ and
		$ \tilde{b}_{k+1},\tilde{b}_{k+2},\ldots,\tilde{b}_{s}=b_{k+1},b_{k+2},\ldots,b_s $. Let $ \tilde{h}>h $ be the case. Clearly, $ h+1<\tilde{h}+1\leq \tilde{k} $. If $ k<\tilde{k} $, we have $ \tilde{b}_{\tilde{k}}=b_{\tilde{k}}<a_{\tilde{h}}=\tilde{a}_{\tilde{h}} $ contradicts  the definition of $ (h,k) $. If  $ h+1<\tilde{h}+1\leq \tilde{k}\leq k  $, from $ b_{k}<a_h<a_{\tilde{h}}=\tilde{a}_{\tilde{h}} $ we have a contradiction with the definition of $ (h,k) $. Therefore, $ \tilde{h}\leq h $. Suppose $ \tilde{h}=h $. If $ \tilde{k}>k $ $ (\tilde{k}=k) $ we have $ b_{\tilde{k}}=\tilde{b}_{\tilde{k}}<\tilde{a}_{\tilde{h}}=\tilde{a}_{h}=b_{k} $ $ (a_h=\tilde{b}_{k}=\tilde{b}_{\tilde{k}}<\tilde{a}_{\tilde{h}}=a_{\tilde{h}}=b_{k}) $ which is a contradiction. Therefore, $ \tilde{k}<k $. This proves that the induction on $ (h,k) $ yields.
	\end{proof}
\end{lem}

\begin{lem}\label{col_reduction}
	A diagonal sorted pair like $ \diagSeqA\geq_c \diagSeqB $ always transforms to a column-wise sorted form like $ \diagSeqC\geq_c \diagSeqD $ of the same shape. 
	\begin{proof}
			If $ \diagSeqA\geq_c \diagSeqB $ is not column-wise sorted, there exists $ 1\leq h \leq \Omega_c $ such that $ a_h > b_h  $ and $ a_h\notin \interval{b} $. Let $ h $ be the maximum index with this property. There exists $ 0\leq v \leq h-1 $  such that $ a_{h-i}>b_{h-i} $  for all  $ 0 \leq i \leq v $ and $ a_{h-v-1}\leq b_{h-v-1} $ if $ v\neq h-1 $. By the definition of $ h $, it is clear that $ a_{h-v-1}\leq b_{h-v-1}<_1 b_{h-v} $ and $ b_h<a_h<_1 a_{h+1} $. Therefore, $ \tilde{a}=a_1,\ldots,a_{h-v-1},b_{h-v},\ldots,b_{h},a_{h+1},\ldots,a_s $ is a chain. 	From the definition of column-wise relations and $ \Omega_c $, one can see that $ (\sigma,\tilde{a}) $ is a well-defined labeled chain. On the other hand, $ b_{h-v-1}<_1 b_{h-v}<a_{h-v} $ and $ a_h<_1 b_{h+1} $ from $ a_h\notin\interval{b} $ and the definition of $ h $. Therefore, $ \tilde{b}=\chain{b}{h-v-1},a_{h-v},\ldots,a_{h},b_{h+1},b_s $ is a chain. From the definition of column-wise relations and $ \Omega_c $, one can deduce that $ (\gamma,\tilde{b}) $  is a well-defined labeled chain. We need to show that $ (\sigma,\tilde{a})\geq_c(\gamma,\tilde{b}) $ is diagonal sorted and the induction on $ h $ converges. \\
			\textbf{Diagonal relations:} Let $ \tilde{a}_i>\tilde{b}_j $ for some $ 1\leq i \leq \Omega_d $ and $ i+1\leq j \leq s $ and $ \tilde{a}\notin\interval{\tilde{b}} $. The case $ h-v\leq i \leq \Omega_d $ leads to trivial contradictions. From the definition of $ v $ in the previous part, it implies $ a_i<b_{h-v-1} $. Hence $ 1\leq i \leq h-v-3 $ and $ i+1\leq j \leq h-v-2 $. So, $ \tilde{a}_i=a_i $ and $ \tilde{b}_j=b_j $. Therefore, from the definition of $ (\sigma,a)\geq_c(\gamma,b) $,  there exists $ j+1 \leq t \leq s $ such that $ b_t-1\leq a_i \leq b_t $. Note that, $ a_i<b_{h-v-1} $ again implies $ t\leq h-v-2 $. Therefore, $ b_t=\tilde{b}_t $. This implies $ a_i=\tilde{a}_i\in\interval{\tilde{b}} $ which is a contradiction. Hence, any diagonal relation in  $ (\sigma,\tilde{a})\geq_c(\gamma,\tilde{b}) $ admits a contradiction. \\
			\textbf{Convergence:} Let $ \tilde{h} $ be the analogues of $ h $  for the new pair $ (\sigma,\tilde{a})\geq_c(\gamma,\tilde{b}) $. From the definition of $ h $, $ \tilde{a} $, $ \tilde{b} $ and the fact that  $ (\sigma,\tilde{a})\geq_c(\gamma,\tilde{b})  $ is diagonal sorted, it is clear that $ \tilde{h}<h $. Therefore, the reverse induction on $ h $ yields the pair $ \diagSeqC \geq_c \diagSeqD $. Since every step of the above process gives a smaller $ h $  and the pair obtained at every step remains diagonal sorted, one can deduce that $ \diagSeqC \geq_c \diagSeqD $ is column-wise sorted. 
	\end{proof}
\end{lem}		

\begin{lem}\label{anti_diagonal_reduction}
	A column-wise sorted pair like $ \diagSeqA\geq_c \diagSeqB $ always transforms to an anti-diagonal sorted (standard form) like $ \diagSeqC\geq_c \diagSeqD $ of the same shape.
	\begin{proof}
		If $ \diagSeqA\geq_c \diagSeqB $ is not anti-diagonal sorted, then there exists   $ 1\leq h \leq \Omega_{ad} $ and $ h+1\leq k \leq r $ such that $ b_h>a_{k} $ and $ b_h\notin\interval{a} $. Assume $ (h,k) $ is minimum with respect $ \leq_\tau $. Recall that $ \leq_\tau $ is $ \operatorname{lex} $ order induced by $ x_1>x_2>\ldots>x_n $. There exists $ 0\leq v \leq h-1 $, such that $ b_{h-i}>a_{k-i}  $ for all $ 0 \leq i \leq v $ and $ b_{h-v-1}\leq a_{k-v-1} $ if $ v\neq h-1 $. We have $ a_{k-v-1}<_1 a_{k-v} <b_{h-v} $. 	From $  b_h\notin\interval{a}  $ and the definition of $ \Omega_{ad} $, we have $ b_h<_1 a_{k+1}  $ (or $ b_h<_1 \sigma_2  $). Therefore, $ \tilde{a}=\chain{a}{k-v-1},b_{h-v},\ldots,b_{h},a_{k+1},\ldots,a_r $  is a chain and  $ (\sigma,\tilde{a}) $ is a well-defined labeled chain. On the other hand, we have $ b_{h-v-1}\leq a_{k-v-1}<_1 a_{k-v} $. Moreover, $ a_{k}<b_h<_1 b_{h+1} $ (or $ a_{k}<b_h\leq \gamma_2 $ when $ h=s $). Thus, $ \tilde{b}=\chain{b}{h-v-1},a_{k-v},\ldots,a_{k},b_{h},\ldots,b_s $ is a chain and $ (\gamma,\tilde{b}) $ is a well-defined labeled chain. To complete the proof, we need to show that $ (\sigma,\tilde{a})\geq_c(\gamma,\tilde{b}) $  is  diagonal sorted and column-wise sorted. Moreover, we need to show that the induction on $ (h,k) $ converge. \\
		\textbf{Diagonal relations:} Let $ \tilde{a}_i> \tilde{b}_j $ with $ 1\leq i \leq \Omega_{d} $ and $ i+1\leq j \leq s $ and $ \tilde{a}_i\notin\interval{\tilde{b}} $. Assume $ (i,j) $ is minimum with respect ot $ \leq_\tau $. Recall that $ v $ is defined in the first part of the proof.  The case $ h-v-1\leq i \leq \Omega_d $ gives trivial contradictions. Suppose $ 1\leq i \leq h-v-2  $. The case $ h-v \leq j \leq s $ also gives trivial contradictions. Suppose $ 1\leq i \leq h-v-2 $ and $ i+1 \leq j \leq h-v-1 $. This implies $ \tilde{a}_i=a_i $  and $ \tilde{b}_j=b_j $. From the definition of $ (\sigma,a)\geq_c(\gamma,b) $ and $ a_i=\tilde{a}_i>\tilde{b}_j=b_j $ we have $ a_i\in\interval{b} $. Hence $ b_t-1\leq a_i \leq b_t $ for unique $ i+1\leq t \leq s $. If $ h<t\leq s  $, we have a contradiction with $ \tilde{a}_i\notin\interval{b} $. If $ h-v\leq t \leq h $, there exists a unique $ 0\leq t'\leq v $ such that $ h-t'=t $ and $ b_{t}>a_{k-t'} $. This implies $ b_t-1\leq a_i<_1 a_{h-t'}<b_t $ which is a contradiction. Hence $ 1\leq t \leq h-v $ which is again a contradiction with $ \tilde{a}_i\notin\interval{\tilde{b}} $. Therefore existence of any diagonal relation in $ (\sigma,\tilde{a})\geq_c(\gamma,\tilde{b}) $ leads to a contradiction.\\
		\textbf{Column-wise relations:} We have already seen that $ (\sigma,\tilde{a})\geq_c(\gamma,\tilde{b}) $ is diagonal sorted. Let $ \tilde{a}_i>\tilde{b}_i $ for $ 1\leq i \leq \Omega_c $ and $ \tilde{a}_i\notin\interval{\tilde{b}} $. Assume $ i $ is maximum. The case $ h-v \leq i \leq \Omega_c $  implies trivial contradictions. Suppose $ 1\leq i \leq h-v-1 $. This gives $ \tilde{a}_i=a_i $ and $ \tilde{b}_i=b_i $. Hence from the definition of $ (\sigma,a)\geq_c(\gamma,b) $, we have $ a_i\in\interval{b} $. Therefore for unique $ i+1\leq t \leq s $, we have $ b_t-1\leq a_i \leq b_t $. The case $ h\leq t \leq s $ implies a contradiction with $ \tilde{a}_i\notin\interval{\tilde{b}} $. Suppose $ h-v \leq t \leq h $. Hence there is a unique $ 0\leq t' \leq v $ such that $ h-t'=t $ and $ a_{k-t'}<b_t $. Therefore, $ b_t-1 \leq a_i <_1  a_{h-t'} < b_t $. This bound  contradicts the definition of $ <_1 $. Hence, $ 1\leq t \leq h-v $ which, again, contradicts $ \tilde{a}_i\notin\interval{\tilde{b}} $. Therefore any column-wise relation in $ (\sigma,\tilde{a})\geq_c(\gamma,\tilde{b}) $ leads to a contradiction. \\
		\textbf{Convergence:} It remains to show that the induction on $ (h,k) $ converges. Let $ (\tilde{h} , \tilde{k}) $ be the analogous of $ (h , k) $ for the new pair $ (\sigma,\tilde{a})\geq_c(\gamma,\tilde{b}) $. It is enough to show $ (\tilde{h} , \tilde{k}) \geq_\tau (h , k) $. By the construction of $ \tilde{a} $ and $ \tilde{b} $, we have $ \tilde{a}_{k+1},\tilde{a}_{k+2},\ldots,\tilde{a}_{r}=a_{k+1},a_{k+2},\ldots,a_{r} $ and $ \tilde{b}_{h+1},\tilde{b}_{h+2},\ldots,\tilde{b}_{s}=b_{h+1},b_{h+2},\ldots,b_{s} $. 	Let $ \tilde{h}>h $ be the case. Clearly, $ h+1<\tilde{h}+1\leq \tilde{k} $. If $ k<\tilde{k} $, we have $ \tilde{a}_{\tilde{k}}=a_{\tilde{k}}<b_{\tilde{h}}=\tilde{b}_{\tilde{h}} $ which is  a contradiction with the definition of $ (h,k) $. If  $ h+1<\tilde{h}+1\leq \tilde{k}\leq k  $, from $ a_{k}<b_h<b_{\tilde{h}} $ we have a contradiction with the definition of $ (h,k) $.  Therefore, $ \tilde{h}\leq h $. Suppose $ \tilde{h}=h $. If $ \tilde{k}>k $ (or $  \tilde{k}=k $) we have $ a_{\tilde{k}}=\tilde{a}_{\tilde{k}}<\tilde{b}_{\tilde{h}}=\tilde{b}_{h}=a_{k} $ (or $ b_h=\tilde{a}_{k}=\tilde{a}_{\tilde{k}}<\tilde{b}_{\tilde{h}}=\tilde{b}_{h}=a_{k} $) which is a contradiction with the definition of $ (h,k) $. Therefore, $ \tilde{k}<k $. This proves that the induction on $ (h,k) $ converges. 
		\end{proof}
\end{lem}

\begin{exmp}
	The tabel
	\[
	\mathcal{A}=\begin{pmatrix}
	(1, 30)       & \vline & 1& 4& 18& 24& 30 &  \\
	(2, 29)       & \vline & 5& 7& 11& 15& 17& 19& 22& 28  \\
	\end{pmatrix}
	\]
	is not standard. We have $ \Delta=1 $ and $ (4,7) $ are the analogues  of $ (h,k) $ is the proof of \autoref{diagonal_reduction}. After applying \autoref{diagonal_reduction}, we obtain the diagonal sorted tabel
	\[
	\mathcal{A}^\prime=\begin{pmatrix}
	(1, 30)       & \vline & 1& 4& 15& 18& 30 &  \\
	(2, 29)       & \vline & 5& 7& 11& 17& 19& 22& 24& 28  \\
	\end{pmatrix}
	\]
	The tabel $ \mathcal{A}' $ has column-wise relations. The analogous of $ h $ in the proof of \autoref{col_reduction} is $ 3 $. By applying \autoref{col_reduction}, we obtain the diagonal sorted and column-wise sorted tabel
	\[
	\mathcal{A}^{\prime\prime}=\begin{pmatrix}
	(1, 30)       & \vline & 1& 4& 11& 18& 30 &  \\
	(2, 29)       & \vline & 5& 7&15& 17& 19& 22& 24& 28  \\
	\end{pmatrix}
	\]
	This tabel has anti-diagonal relations. The analogous of $ (h,k) $ in the proof of \autoref{anti_diagonal_reduction} is $ (1,2) $. Finally, by applying \autoref{anti_diagonal_reduction}, we obtain the following standard tabel:
	\[
	\mathcal{B}=\begin{pmatrix}
	(1, 30)       & \vline & 1& 5& 11& 18& 30 &  \\
	(2, 29)       & \vline & 4& 7& 15& 17& 19& 22& 24& 28  \\
	\end{pmatrix}
	\]
\end{exmp}

Let $ \tabelA{\sigma}{a} $ be a tabel. The coordinates of the "cells" are denoted by $ (i,j) $ where $ i $ and $ j $ are row and column indices respectively.

\begin{algorithm}\label{label_algorithm}
	Let  $ \tabelA{\sigma}{a} $ be a tabel. We assign a label to each coordinate of $ \mathcal{A} $ by starting from the cell with coordinate $ (1,1) $ (i.e first row and first column) and increasing the labels as we go through every cell in the first column. Then we proceed by increasing the labels for next columns.  Let $ \operatorname{P}(i,j) $ be the function that assigns a label to the coordinate $ (i,j) $ where $ i $ is the row index and $ j $ is the column index. 
\end{algorithm}

\begin{exmp}
	Let $ \mathcal{A} $ be a tabel with shape $ (\sigma_1,8),(\sigma_2,4),(\sigma_3,6),(\sigma_4,5) $. The \autoref{label_algorithm}, labels the cells of $ \mathcal{A} $ like the following:
	$$ \mathcal{A}= \left(\begin{array}{*{20}c}
	\sigma_1 & \vline &  1 & 5 & 9 & 13 & 17 & 20 & 22 & 23 &\\
	\sigma_2 & \vline & 2 & 6 & 10 & 14 &  &  & &\\
	\sigma_3 & \vline & 3 & 7 & 11 & 15 & 18 & 21 &\\
	\sigma_4 & \vline & 4 & 8 & 12 & 16 & 19 &  &  &  & &
	\end{array}
	\right)$$	
\end{exmp}

 We say the entry $ a^{(i)}_j $ in $ \mathcal{A} $ is stable modulo diagonal relations if the coordinate of $ a^{(i)}_j $ remains the same after any pairwise diagonal relations transformations of rows of $ \mathcal{A} $. Let stable entries modulo column-wise relations and anti-diagonal relations be defined accordingly. It is clear that $ \mathcal{A} $ is standard if and only if all the entries of $ \mathcal{A} $ are stable modulo diagonal relations, column-wise relations and anti-diagonal relations.

\begin{prop}\label{tabel_reduction}
	Let $ \tabelA{\sigma}{a} $ be a tabel. Then, $ \mathcal{A} $ always reduces to a standard form with the same shape. 
	\begin{proof}
		We treat transformation modulo each type of relations separately. Let $ c_1\leq \ldots \leq c_l $ be the entries of $ \mathcal{A} $ in order where the entries are labeled by \autoref{label_algorithm}. We argue by revers induction on $ 1\leq t \leq l $.  Suppose $ \mathcal{A}  $ has some diagonal relations. 
		\begin{itemize}
			\item Let $ t=l $. Note that since $ c_l $ is the largest entry and the rows of $ 
			\mathcal{A} $ are chains, this entry is located at the last cell of some row. Let $ (i,r_i) $ be the coordinate of  $ c_l $.  Suppose $ c_l=a^{(i)}_{r_i} $ is not stable modulo diagonal relations.   There exists some row index $ i<j $ such that $ a^{(i)}_{r_i}>a^{(j)}_{r_j} $ and $ a^{(i)}_{r_i} \leq \sigma_{j,2} $ and of course $ r_i<r_j $. Let $ j $ be maximum. By applying \autoref{diagonal_reduction} on the rows $ (i,j) $, one can reduce the tabel into a new tabel in which $ c_l $ is stable modulo diagonal relations by construction. 
			\item Let $ c_{t+1} $ be stable modulo diagonal relations for $ t<l $. Let $ (i,h) $ be the coordinate of $ c_t=a^{(i)}_{h} $. From the definition of diagonal relations, there exists some row index $ i<j $ such that $ a^{(i)}_{h}>a^{(j)}_{k} $ and $ a^{(i)}_{h}\notin\interval{a^{(j)}} $ for some $ k $ with $ h+1\leq k \leq \Omega_{d}(i,j) $. Assume $ j $ is maximum. By applying \autoref{diagonal_reduction} on the rows $ (i,j) $, one can reduce to a new tabel in which $ c_t $ is stable modulo diagonal relations by construction. 
		\end{itemize}
	Let $ \mathcal{A} $ be obtained form the previous step. Suppose $ \mathcal{A} $ have some column-wise relations. 
	\begin{itemize}
		\item Let $ t=l $. Note that since $ c_l $ is the largest entry and the rows of $ 
		\mathcal{A} $ are chains, this entry is located at the last cell of some row. 	Let $ (i,r_i) $ be the coordinate of $ c_l=a^{(i)}_{r_i} $.  There exists a row index $ j $ such that $ r_i=r_j $, $ a^{(i)}_{r_i}>a^{(j)}_{r_i} $ and $ a^{(i)}_{r_i}\leq \sigma_{j,2} $. Let $ j $  be  maximum. By applying \autoref{col_reduction} on the rows $ (i,j) $ and replacing it in $ \mathcal{A} $, one can reduce to a new tabel in which $ c_l $ is stable modulo diagonal relations and column-wise relations by construction.
		\item Let $ c_{t+1} $ be stable for $ t<l $. Suppose $ c_t $ is not stable modulo column-wise relations. Let $ (i,h) $ be the coordinate of $ c_t=a^{(i)}_{h} $.
		There exists some row $ j $ such that $ a^{(i)}_{h}>a^{(j)}_{h} $ and $ a^{(i)}_{h}\notin\interval{a^{(j)}} $. Assume that $ j $ is maximum. By applying \autoref{col_reduction} on the rows $ (i,j) $ and replacing it in $ \mathcal{A} $, one can reduce to a tabel in which $ c_t $ is  stable modulo diagonal relations and column-wise relations by construction.
	\end{itemize}
		Finally, assume $ \mathcal{A} $ is obtained by applying last two steps. Meaning that $ \mathcal{A} $ is diagonal sorted and column-wise sorted. Suppose $ \mathcal{A} $ has anti-diagonal relations. 
		\begin{itemize}
			\item Let $ t=l $. Note that since $ c_l $ is the largest entry and the rows of $ 
			\mathcal{A} $ are chains, this entry is located at the last cell of some row. 	Let $ (i,r_i) $ be the coordinate of $ c_l=a^{(i)}_{r_i} $.  There exists a row index $ j<i $ such that $ a^{(i)}_{r_i}>a^{(j)}_{r_j} $ and $ a^{(i)}_{r_i}\leq \sigma_{j,2} $ and of course $ r_i<r_j $. Let $ r_j $ be maximum and  $ j $ be maximum given $ r_j $.  By applying \autoref{anti_diagonal_reduction} on the rows $ (i,j) $ and replacing it in $ \mathcal{A} $, one can reduce to a new tabel in which  $ c_l $ is stable modulo diagonal, column-wise and anti-diagonal relations by construction.
			\item Let $ c_{t+1} $ be stable for $ t<l $. Let $ (i,h) $ be the coordinate of  $ c_t=a^{(i)}_{h} $.  There exists some row index $ j<i $ such that $ a^{(i)}_{h}>a^{(j)}_{k} $ for some $ h+1\leq k \leq \Omega_{ad}(i,j) $ and $ a^{(i)}_{h}\notin\interval{a^{(j)}} $. Let $ k $ be maximum and $ j $ be maximum given $ k $.  By applying \autoref{anti_diagonal_reduction} on the rows $ (i,j) $ and replacing it in $ \mathcal{A} $, one can reduce to a new tabel in which  $ c_t $ is stable modulo diagonal, column-wise and anti-diagonal relations by construction.
		\end{itemize}
	\end{proof}
\end{prop}

\begin{lem}\label{reduction_uniquness}
	Let $ \tabelA{\sigma}{a} $ be a tabel. Then, $ \mathcal{A} $ always reduces to a unique standard form $ \tabelB{\sigma}{b} $  of the same shape. 
	\begin{proof}
		Let $ \mathcal{A} $ be labeled by \autoref{label_algorithm}. We have $ \sum_{i}r_{i}=l  $. Let $ \mathcal{B} $ be a standard form reduction of $ \mathcal{A} $. Since $ \mathcal{A} $ and $ \mathcal{B} $ have the same shape, the \autoref{label_algorithm} assigns the same labels to $ \mathcal{B} $.  Let $ c_1\leq c_2\leq \ldots \leq c_l $ be the ordered set of the entries of $ \mathcal{A} $ and $ \mathcal{B} $.  Note that the function in \autoref{label_algorithm} is a one to one correspondence. We prove that for all $ t\colon 1, \ldots, l $, the entries $ b^{(i)}_j $ with $ \operatorname{P}(i,j)=t $ is determined uniquely. Therefore, $ \mathcal{B} $ is given uniquely. We proceed by revers induction on $ t $.

		\begin{enumerate}[(1)]
			\item Let $ t=l $. There exists unique coordinate $ (i,j) $ such that $ \operatorname{P}(i,j)=l $. Recall that by \autoref{label_algorithm}, $ (i,j) $ is the coordinate of the last cell of the longest row.  If $ c_l\leq \sigma_{i,2} $, from \autoref{label_algorithm} and part $ (2) $ of the \autoref{largest_row_last_entry}, we have $ b^{(i)}_j=c_l $. If $ c_l> \sigma_{i,2} $, from the fact that $ \mathcal{B} $ is a well-defined tabel, there exists $ c_h $ such that $ c_h<c_{h+1}=\ldots=c_l $. Thus $ b^{(i)}_j\leq c_h \leq \sigma_{i,2} $. If $ b^{(i)}_j<c_h $, there exists some row index of tabel $ \mathcal{B} $ like $ i^\prime\neq i $ which contains $ c_h $. This yields that there exists either diagonal relations, column-wise relations or anti diagonal relations in $ \mathcal{B} $. Which contradicts the definition of $ \mathcal{B} $.  Hence, $ b^{(i)}_j=c_h $.
			\item Let uniqueness of  $ b^{(i')}_{j'} $ be given for every coordinate with $ t<\operatorname{P}(i^\prime,j^\prime)\leq l $. Let $ c_1\leq c_2\leq \ldots \leq c_t $ be the remainder of the entries of $ \mathcal{A} $ and $ \mathcal{B} $ relabeled by $ 1,2,\ldots ,t $. There exists a unique coordinate $ (i,j) $  with $ \operatorname{P}(i,j)=t $. We always have $ b^{(i)}_j\leq c_t $ .
			\begin{itemize}
				\item If $ j=r_i $ and $ c_t\leq\sigma_{i,2} $.  There exists some coordinate $ (i_t,j_t) $ with $ \operatorname{P}(i_t,j_t)\leq t $ and $ b^{(i_t)}_{j_t} = c_t $.  If $ b^{(i)}_j< c_t $, we have $ i_t\neq i $ since $ b^{(i)} $ is a chain. According to the induction hypothesis and the \autoref{label_algorithm}, we have $ i_t< i $ and $ j_t\leq j $ or $ i_t>i $  and $ j_t>j $. This means that $ (\sigma_{i_t},b^{(i_t)}) $ and $ (\sigma_i,b^{(i)}) $ have either diagonal relations, column-wise relations or anti diagonal relations. This is a contradiction with the definition of $ \mathcal{B} $. So $ b^{(i)}_j= c_t $. \\
				If $ c_t>\sigma_{i,2} $, since $ \mathcal{B} $ is a well-defined tabel, we can find $ c_{t^\prime} $ where  $ c_{t^\prime}<c_{t^\prime+1}=\ldots=c_t $. We have $b^{(i)}_j\leq c_{t^\prime} $. Let $ b^{(i)}_j< c_{t^\prime} $. There exist a unique coordinate $ (i_{t^\prime},j_{t^\prime}) $ with $ \operatorname{P}(i_{t^\prime},j_{t^\prime})\leq t $ such that $ b^{(i_{t^\prime})}_{j_{t^\prime}}=c_{t^{\prime}} $. According to the induction hypothesis and the \autoref{label_algorithm}, we have $ i_{t^\prime}< i $ and $ j_{t^\prime}\leq j $ or $ i_{t^\prime}>i $  and $ j_{t^\prime}>j $. This means that $ (\sigma_{i_{t^\prime}},b^{(i_{t^\prime})}) $ and $ (\sigma_i,b^{(i)}) $ have either diagonal relations, column-wise relations or anti diagonal relations. This is a contradiction with the definition of $ \mathcal{B} $. So $ b^{(i)}_j= c_{t^\prime} $.
				\item If $ j\neq r_i $ and $ c_t+1< b^{(i)}_{j+1} $. There exists a unique coordinate $ (i_t,j_t) $ with $ \operatorname{P}(i_t,j_t)\leq t $ and $ b^{(i_t)}_{j_t}=c_t $. If $ b^{(i)}_j<c_t $, we have $ i_t\neq i $ since $ b^{(i)} $ is a chain. According to the induction hypothesis and the \autoref{label_algorithm}, we have $ i_t< i $ and $ j_t\leq j $ or $ i_t>i $  and $ j_t>j $. From $ c_t+1< b^{(i)}_{j+1} $ and $ b^{(i)}_j<c_t $, we have $ c_t\notin\interval{b^{(i)}} $. This means that $ (\sigma_{i_{t}},b^{(i_{t})}) $ and $ (\sigma_i,b^{(i)}) $ have either diagonal relations, column-wise relations or anti diagonal relations. This is a contradiction with the definition of $ \mathcal{B} $. So $ b^{(i)}_j= c_{t} $.\\	
				If $ c_t+1\geq b^{(i)}_{j+1} $, from the fact that $ \mathcal{B} $ is a well-defined tabel, we can find $ c_{t^\prime}<_1b^{(i)}_{j+1} $. Let $ t^{\prime} $ be the largest label satisfying this condition. Hence, $ b^{(i)}_j\leq c_{t^\prime} $. There exists a unique coordinate $ (i_{t^\prime},j_{t^\prime}) $ with $ \operatorname{P}(i_{t^\prime},j_{t^\prime})\leq t $ such that $ b^{(i_{t^\prime})}_{j_{t^\prime}}=c_{t^\prime} $. If $ b^{(i)}_j< c_{t^\prime} $, we have $ i_{t^\prime}\neq i $ since $ b^{(i)} $ is a chain. According to the induction hypothesis and the \autoref{label_algorithm}, we have $ i_{t^\prime}< i $ and $ j_{t^\prime}\leq j $ or $ i_{t^\prime}>i $  and $ j_{t^\prime}>j $. Since $ b^{(i)}_j<c_{t^\prime}<_1b^{(i)}_{j+1} $, we have $ c_{t^\prime}\notin\interval{b^{(i)}} $. This means that $ (\sigma_{i_{t^\prime}},b^{(i_{t^\prime})}) $ and $ (\sigma_i,b^{(i)}) $ have either diagonal relations, column-wise relations or anti diagonal relations. This is a contradiction with the definition of $ \mathcal{B} $. So $ b^{(i)}_j= c_{t^\prime} $.

			\end{itemize}	
		\end{enumerate}		
	\end{proof}
\end{lem}

\begin{exmp}
	The  non-standard tabel
	$$ \mathcal{A}= \left(\begin{array}{*{20}c}
	(1, 29) & \vline &  8 & 12 & 18 & 20 & 22 &  \\
	(1, 30) & \vline & 2 & 7 & 23 & 25 & 27 & 30 & \\
	(1, 30) & \vline & 1 & 18 & 23 & 27 & 30 &\\
	(2, 29) & \vline & 2 & 5 & 7 & 9 & 13 & 16 & 20 & 25 & \\
	(2, 30) & \vline & 8 & 10 & 12 & 17 & 25 & 28 & 
	\end{array}
	\right)$$
	transforms to the following standard tabel: 
	$$ \mathcal{B}= \left(\begin{array}{*{20}c}
	(1, 29) & \vline & 1 & 7 & 12 & 18 & 23 &  \\
	(1, 30) & \vline & 2 & 8 & 12 & 20 & 25 & 30 & \\
	(1, 30) & \vline & 2 & 8 & 13 & 20 & 27 &\\
	(2, 29) & \vline & 5 & 9 & 16 & 18 & 20 & 23 & 25& 28 & \\
	(2, 30) & \vline & 7 & 10 & 17 & 22 & 27 & 30 & 
	\end{array}
	\right)$$
\end{exmp}

\begin{rem}\label{ugly_remark}
	As we saw, the function $ \Delta $ in fact controls the well-definity of the transformations of labeled chains with respect to our relations. Moreover, $ \geq_c $ decides the order of the rows of the tabels. In other words, for a given pair of chains $ a=\chain{a}{r} $ and $ b=\chain{b}{s} $, we can make a tabel with first row $ a=\chain{a}{r} $ and second row $ b=\chain{b}{s} $ with out considering any labels. Let us  denote this tabel by $ (a,b) $. Now, by omitting the role of  $ \Delta $  by setting $ \Delta=0 $, we can always perform \autoref{diagonal_reduction},\autoref{col_reduction}, \autoref{anti_diagonal_reduction} and \autoref{reduction_uniquness}. In particular, the following holds:
	\begin{enumerate}[(I)]
		\item The tabel $ (a,b) $ is standard if and only if 
		\begin{enumerate}[(i)]
			\item  $ a_i\leq b_{i+1} $ for $ i\colon 1,\ldots,\Omega_{\text{d}} $ or  $ a_h>b_{h+1} $ for some $ h $ and $ a_h,\ldots,a_{\Omega_{\text{d}}}\in \interval{b} $,
			
			\item  $ a_i\leq b_{i} $ for $ i\colon 1,\ldots,\Omega_{\text{c}} $ or  $ a_h>b_{h} $ for some $ h $ and $ a_h,\ldots,a_{\Omega_{\text{c}}}\in \interval{b} $,
			
			\item   $ b_{i}\leq a_{i+1}  $ for $ i\colon 1,\ldots,\Omega_{\text{ad}} $ or  $ b_h>a_{h+1} $ for some $ h $ and $ b_h,\ldots,b_{\Omega_{\text{ad}}}\in \interval{a} $.
		\end{enumerate}
		\item In particular, when $ r\geq s $, the tabel $ (a,b) $ is standard if and only if
		\begin{enumerate}[(i)]
			\item $ a_i\leq b_{i} $ for all $ 1\leq i \leq s $;
			\item $ b_i\leq a_{i+1} $ or $ b_h>a_{h+1} $ for some $ h $ and $ b_h,\ldots,b_{s}\in\interval{a}. $\\
		\end{enumerate}
		\item If  $ (a,b)  $ is standard, we have $ a_r\leq b_s $ when $ r\leq s $ and $ b_s\leq a_r $ when $ r> s $. 
	\end{enumerate}
\end{rem}

\section{Sagbi Deformations and Multi-Rees Algebra}\label{Multi-Rees_Algebra}
	In this section, we use the machinery introduced in  Section \ref{Standard_Forms} to study the multi-Rees algebra of ideals of  family $ \f=\f^{(1,n)}\cup\f^{(1,n-1)}\cup\f^{(2,n)}\cup\f^{(2,n-1)} $. Let $ \chain{I}{l} $ be ideals of the ring $ S=K[\bar{x}] $ where $ \bar{x}:=\chain{x}{n} $. The multi-Rees algebra of $ \chain{I}{l} $ is 
	$$ \rees{\chain{I}{l}}=\bigoplus_{\chain{\alpha}{l}}(I_1t_1)^{\alpha_1}\ldots(I_lt_l)^{\alpha_l} $$
	where  $ \bar{t}:=\chain{t}{l} $ are new indeterminates over $ S $ and $ (\chain{\alpha}{l})\in\mathbb{N}^l $. One can also see multi-Rees algebra of $ \chain{I}{l} $ as 
	$$ \rees{\chain{I}{l}}=S[I_1t_1,\ldots,I_lt_l]\subset S[\,\bar{t}\,] $$
	\begin{nota}
		In this work, we denote a maximal minor of the Hankel matrix $ \hmatTagSize{\sigma}{r} $ by $ [\chain{a}{r}] $ where $ \chain{a}{r} $ is the chain of the entries on the main diagonal. Moreover, we use $ \detIdeal{\sigma}{r} $  to denote the determinantal ideal of maximal minors of $ \hmatTagSize{\sigma}{r} $. 
	\end{nota}
	Let $ \bar{t} $ be the set of all new indeterminates $ t_{\sigma,r} $ over $ S $ where $ \sigma $ and $ r $ go through  labels and  lengths of all labeled chains $ (\sigma,a) $ with $ a=\chain{a}{r} $. Consider the multi-Rees algebra $ \reesSymbol=\rees{I_{\sigma,r}t_{\sigma,r} \colon I_{\sigma,r}\in\f} $ of all ideals of the family $ \f $ and the multi-Rees algebra $ \reesSymbol^{\operatorname{in}}=\rees{\initial{I_{\sigma,r}}t_{\sigma,r} \colon I_{\sigma,r}\in\f} $. We shall consider $ \reesSymbol $ and $ \reesSymbol^{\operatorname{in}} $ as sub rings of $ S[\,\bar{t}\,] $. One can also consider the representation of our Rees algebras as quotients of some polynomial ring.\\   
	Let $ \bar{z} $ be the set  of new indeterminates $ z_{\sigma,a} $ over $ S $  where $ (\sigma,a) $ runs through all labeled chains. 	Consider the polynomial ring $ R=S[\,\bar{z}\,] $.  Recall that we reserve the letters $ r $ and $ s $ for the lengths of chains $ a=\chain{a}{r} $ and $ b=\chain{b}{s} $ respectively. 
	Consider the following surjective algebraic homomorphisms:
	\begin{align*}
		\varphi\colon R\rightarrow& \rees{I_{\sigma,r}t_{\sigma,r} \colon I_{\sigma,r}\in\f}\\
		x_i \mapsto & x_i \\
	z_{\sigma,a} \mapsto & [a]t_{\sigma,r} \\
	\end{align*}
	and 
	\begin{align*}
		\varphi^{\operatorname{in}}\colon R\rightarrow& \rees{\initial{I_{\sigma,r}}t_{\sigma,r} \colon I_{\sigma,r}\in\f}\\
		x_i \mapsto & x_i \\
		z_{\sigma,a}  \mapsto & x_{a}t_{\sigma,r}
	\end{align*}
	
It is known that the maximal minors of $ \hmatTagSize{\sigma}{r} $ form a Gr\"{o}bner basis for  ideal $ I_{\sigma,r} $ with respect to any diagonal term order. Therefore $ \varphi^{\operatorname{in}} $ is surjective. \\
Often the structure of multi-Rees algebras are better understood by looking at their representation as a quotient of a polynomial ring. Consider isomorphisms $ R/\kerRees\simeq\reesSymbol $ and $ R/\kerReesIn\simeq\reesSymbol^{\operatorname{in}} $ induced by $ \varphi $ and $ \varphi^{\operatorname{in}} $.\\
Let $ l $ be the cardinality of $ \bar{z} $. We equip $ R $ with  $ \mathbb{Z}\oplus\mathbb{Z}^{l} $  graded setting by considering   $ \operatorname{deg}(x_i)=e $  and $ \operatorname{deg}(z_{\sigma,a})=e_{\sigma,r} $. Note that $ e $ and $ e_{\sigma,r}$'s are the standard basis for $ \mathbb{Z}\oplus\mathbb{Z}^{l} $.  In order to define $ \varphi $ and $ \varphi^{\operatorname{in}} $ as multi-homogeneous algebraic homomorphisms, in  $ S[\,\bar{t}\,] $ we set $ \deg(x_i)=e $ and $ \deg(t_{\sigma,r})=-re+e_{\sigma,r} $. This will set $ \reesSymbol $ and $ \reesSymbol^{\operatorname{in}} $ as standard multi-graded algebras. The multi-graded setting is effective in the proof of the following proposition. 
\begin{prop}\label{construction_standard_tabel_monomial}
	Let $ mv\in R $  be a monomial where $ m $ is a monomial in $ \bar{x} $ and $ v $ is a monomial in $ \bar{z} $. There exists a unique representation $ \varphi^{\operatorname{in}}(mv)=u\prod_{\mathcal{A}}(x_at_{\sigma,r}) $ where $ u $ is a monomial in $ \bar{x} $ and $ \mathcal{A} $ is a tabel of shape $ \deg(v) $
	such that:
	\begin{enumerate}[(i)]
		\item $ \mathcal{A} $ is a standard tabel;
		\item for every indeterminate $ x_i $ in $ u $ and every row $ (\sigma,a) $ in $ \mathcal{A} $, we have:
		\begin{enumerate}[(a)]
			\item  $ i\leq a_1 $ or
			\item  $ i>a_1 $ and either $ \sigma_2<i $ or $ i \in\interval{a} $.
		\end{enumerate}
	\end{enumerate}
	\begin{proof}
		From the definition of $ \varphi^{\operatorname{in}} $ and the multi-graded setting of $ R $ there exists a representation 
		$ \varphi^{\operatorname{in}}(mv)=c\prod_{\mathcal{B}}(x_bt_{\sigma,r}) $ such that $ \mathcal{B} $ has shape $ \deg(v) $ and $ c $ is a monomial in $ \bar{x} $. Assume $ u $ is the maximum of such $ c $'s with respect to $ \leq_\tau $. By virtue of \autoref{reduction_uniquness}, we take $ \mathcal{A} $ to be the unique standard form of $ \mathcal{B} $'s. It remains to show that $ \varphi^{\operatorname{in}}(mv)=u\prod_{\mathcal{A}}(x_at_{\sigma,r}) $ satisfies $ (i) $ and $ (ii) $. The condition $ (i) $ is clearly satisfied by the construction of  $ \mathcal{A} $. Let $ x_i $ be an indeterminate of $ u $ and  $ (\sigma,a) $ a row in $ \mathcal{A} $ not satisfying $ (ii) $. Then,  for some $ 1\leq t \leq r $  we have $ a_t<i\leq a_{t+1} $  ( or $ a_r <i \leq \sigma_2 $). Replace $ x_i $ and $ x_{a_t} $ (or $ x_i $ and $ x_{a_r} $) and denote the new tabel with $ \mathcal{A}' $. Consider the presentation 	$ \varphi^{\operatorname{in}}(mv)=x_{a_t}u/x_i\prod_{\mathcal{A}'}(x_{a'}t_{\sigma,r}) $. By virtue of \autoref{reduction_uniquness}, we can assume $ \mathcal{A}' $ is standard. We have $ x_{a_t}u/x_i\geq_{\tau}u $ which is a contradiction with the definition of $ u $. Hence, $ \varphi^{\operatorname{in}}(mv)=u\prod_{\mathcal{A}}(x_{a}t_{\sigma,r}) $ satisfies condition $ (ii) $. 
	\end{proof}
\end{prop}
We will refer to the tabel $ \mathcal{A} $ of the above construction by the standard tabel of $ mv $. 
\begin{dfn}
		Let $ mv\in R $  be a monomial where $ m $ is a monomial in $ \bar{x} $ and $ v $ is a monomial in $ \bar{z} $. We define $ u\prod_{\mathcal{A}}z_{\sigma,a} $ to be the standard form of $ mv $, where $ \mathcal{A} $ and $ u $ are obtained in \autoref{construction_standard_tabel_monomial}. 
	We say $ mv $ is a standard monomial if and only if $ mv = u\prod_{\mathcal{A}}z_{\sigma,a} $. 
\end{dfn}
\begin{rem}\label{rem_on_std_mon}
	The following holds:
	\begin{enumerate}[(i)]
		\item A monomial $ u\prod_{\mathcal{A}}z_{\sigma,a} $ in $ R $ is standard if and only if every factor $ x_iz_{\sigma,a} $ and $ z_{\sigma,a}z_{\gamma,b} $ in $ u\prod_{\mathcal{A}}z_{\sigma,a} $ is standard. 
		\begin{proof}
			Let	$ u\prod_{\mathcal{A}}z_{\sigma,a} $ in $ R $ be standard and let $ x_iz_{\sigma,a} $ and $ z_{\sigma,a}z_{\gamma,b} $ be factors in $ u\prod_{\mathcal{A}}z_{\sigma,a} $. Note that  $ z_{\sigma,a}z_{\gamma,b} $ defines a pair of rows in $ \mathcal{A} $. From \autoref{largest_row_last_entry}, $ z_{\sigma,a}z_{\gamma,b} $  is standard provided $ u\prod_{\mathcal{A}}z_{\sigma,a} $ satisfies \autoref{construction_standard_tabel_monomial} part $ (i) $. Moreover, \autoref{construction_standard_tabel_monomial} part $ (ii) $ clearly stats that $ x_iz_{\sigma,a} $ is standard. Conversely, let every $ x_iz_{\sigma,a} $ and $ z_{\sigma,a}z_{\gamma,b} $ be standard monomials. From \autoref{largest_row_last_entry}, the tabel $ \mathcal{A} $ is standard given every $ z_{\sigma,a}z_{\gamma,b} $ is standard. Now, it remains to show that $ u $ and $ \mathcal{A} $ satisfies \autoref{construction_standard_tabel_monomial} part $ (ii) $. Let that not be the case. There exists  $ x_i $ dividing $ u $ and a row $ (\sigma,a) $ in $ \mathcal{A} $ such that $ a_t<i\leq_c a_{t+1} $ (or $ a_r<i\leq \sigma_2 $). This in fact means that $ x_iz_{\sigma,a} $ is non-standard which contradicts our hypothesis. 
		\end{proof}	
		\item 	Let $ mv $ be a monomial in $ R $ with standard form  $ u\prod_{\mathcal{A}}z_{\sigma,a} $. From  \autoref{construction_standard_tabel_monomial}, \autoref{reduction_uniquness} and the isomorphism $ R/\kerReesIn\backsimeq\reesSymbol^{\operatorname{in}} $, it yields that in a class $ \overline{mv}\in R/\kerReesIn $, there exists exactly one standard monomial. 
	\end{enumerate}
\end{rem}
Consider a marked polynomial to be a polynomial $ f\in R\setminus\{0\} $ together with a specific term $ \initial{f} $ in $ f $. Note that $ \initial{f} $ can be any term of $ f $. For a given finite set of marked polynomials like $ \mathcal{F} $, we define the reduction algorithm modulo $ \mathcal{F} $  in the natural sense. We say that $ \mathcal{F} $ is marked coherently if there exists a term order $ \prec $ on $ R $ such that $ \initial{f}=\operatorname{in}_\prec(f) $ for all $ f \in \mathcal{F} $. It is clear that if $ \mathcal{F} $  is marked coherently, then the reduction modulo $ \mathcal{F} $ is Noetherian. The following is a classic result. 
\begin{theorem}\label{term_order_sturmfels}
	A finite set $ \mathcal{F}\subset R $ of marked polynomials is marked coherently if and only if the reduction modulo $ \mathcal{F} $  converges.  
	\begin{proof}
		\cite[Theorem 3.12]{Sturmfels_marked_polynomials}
	\end{proof}
\end{theorem}  
Consider the following finite set of marked polynomials where the marked terms are underlined.  
\begin{align*}
G=\Big\{&\underline{z_{\sigma,a} z_{\gamma,b}}-z_{\sigma,c} z_{\gamma,d}\; \colon \;  z_{\sigma,a} z_{\gamma,b} \text{ is a non-standard monomial and its standard form is } z_{\sigma,c} z_{\gamma,d}\\
&\underline{x_i z_{\sigma,a}}-x_j z_{\sigma,c}\; \colon \;  x_i z_{\sigma,a} \text{ is a non-standard monomial and its standard form is } x_j z_{\sigma,c} \Big\}
\end{align*}

From \autoref{reduction_uniquness}, \autoref{rem_on_std_mon} and \autoref{term_order_sturmfels}, we can see that there exists a term order on $ R $ which picks the underlined monomials of $ G $ as the leading terms. We denote this term order by $ \leq_\alpha $. 
The following is a classic result:

\begin{lem}\label{G_basis_lem}
	Let $ K[Y_1,\ldots,Y_n] $ be a polynomial ring equipped with some term order. Let $ J\subset K[Y_1,\ldots,Y_n] $  be an ideal and let $ f_1,\ldots,f_s $ be polynomials in $ J $. If the set $ \Omega=\{Y^a \; \colon \; Y^a\notin(\initial{f_1},\ldots,\initial{f_s}) \} $ are linearly independent in $ K[Y_1,\ldots,Y_n]/J $, then $ f_1,\ldots,f_s $  is a Gr\"{o}bner basis of $ J $ with respect to the term order. 
\end{lem}
Now we have everything we need to prove first main theorem of this section. 
\begin{theorem}\label{main_theorem_initial}
	The family $ \f^{\initiaL}=\{\initial{\detIdeal{\sigma}{a}}:\detIdeal{\sigma}{a} \in\f\} $ has the following features:
	\begin{enumerate}[(1)]
		\item Every product of ideals of $ \f^{\initiaL} $ has linear resolution.
		\item The multi-Rees algebra $ \reesSymbol^{\initiaL}=\rees{\initial{I_{\sigma,r}}\colon I_{\sigma,r}\in\f} $ is defined by a quadratic Gr\"{o}bner basis with respect to $ \leq_{\alpha} $  and it is Koszul.
	\end{enumerate}
	\begin{proof}
		It is enough to prove $ (2) $. It is clear that $ G $ is in $ \kerReesIn $. Let $ \Omega=\{mv\colon mv \notin(\initial{g}\colon g\in G)\} $. From \autoref{G_basis_lem}, it is enough to prove that the elements of $ \Omega $ are linearly independent in $ R/\kerReesIn $. Let $ \sum_i\lambda_i \overline{m_iv_i}=0 $ in $ R/\kerReesIn $ where $ m_iv_i\in\Omega $.  \autoref{rem_on_std_mon} part $ (i) $ shows that $ m_iv_i $ is the standard monomial representative of its class.  From   $ R/\kerReesIn\simeq\reesSymbol^{\operatorname{in}}\subseteq S[\,\bar{t}\,] $, we see that $ \overline{m_iv_i} $'s are linearly independent if and only if $ \varphi^{\operatorname{in}}(m_iv_i) $'s are pairwise distinct. This is in fact the case from \autoref{rem_on_std_mon} part $ (ii) $. Hence, $ \lambda_i=0 $ for every $ i $. Thus $ \reesSymbol^{\initiaL} $ is defined by a quadratic Gr\"{o}bner basis. Hence it is Koszul. Now the multi-graded version of the theorem of Blum \cite{Blum} proves $ (1) $. 
	\end{proof}
\end{theorem}
	
In the rest of the section, we apply the means of Gr\"{o}bner basis and sagbi basis to study $ \reesSymbol $. In Section \ref{Standard_Forms}, we saw that the "data" encoded in the product of some labeled chains can be presented as a tabel. We  employ  this tools to lift $ G $ to a Gr\"{o}bner basis for $ \kerRees $. Our main tool to perform the lifting is as simple as the observation of Laplace expansion of the minors. 

\begin{cor}\label{G-basis_pruduct}

	Let $ I_1,\ldots,I_l $  be  ideals of  the family $ \f^{(1,n)} $. Then, the natural generators of $ I=I_1\ldots I_l $ form a Gr\"{o}bner basis with respect to  $ \leq_{\tau} $. In particular for $ l=2 $, if $ \sum_i \lambda_i[a^{(i)}][b^{(i)}] $ is some linear combination, such that $ [a^{(i)}]  $ and $ [b^{(i)}] $ are maximal minors of the matrices $ \hmatTagSize{\sigma}{r} $ and $ \hmatTagSize{\gamma}{s} $ and $ \lambda_i\in K $, then there exist chains $ e=\chain{e}{r} $ and $ f=\chain{f}{s} $ such that 
	$$ \lambda x_ex_f=\initial{\sum_i\lambda_i[a^{(i)}][b^{(i)}]} $$ 
	where $ \lambda\in K $.
	\begin{proof}
		The first part is proved in \cite[Corollary 3.26]{NAM_PAPER}. For the second part, it is enough to consider $ [a^{(i)}]  $ and $ [b^{(i)}] $ as maximal minors of the family $ \hmatTagSize{(1,n)}{} $. Note that this does not affect the polynomials given by this pair of minors. Now from the first part of the statement, the existence of $ e=\chain{e}{r} $ and $ f=\chain{f}{s} $ follows.
	\end{proof}
\end{cor}

\begin{nota}
	Let $ u $ be a monomial   in $ S $. We use $ \deg_{x_i}(u) $ to denote the degree of $ x_i $ in $ u $ (i.e the number of copies of $ x_i $ in $ u $).
\end{nota}

\begin{observe}\label{monomial_x_1_x_n_deg}
	Let $ [a] $ and $ [b] $ be maximal minors of the matrices $ \hmatTagSize{\sigma}{r} $ and $ \hmatTagSize{\gamma}{s} $. Let $ c_1\leq \ldots \leq c_{r+s} $ be the entries of the chains $ a $ and $ b $ in order. Let $ u $ be any term in $ [a][b] $. It is easy to see that $ \deg_{x_{c_1}}(u)\leq \deg_{x_{c_1}}(x_ax_b) $ and $ \deg_{x_{c_{r+s}}}(u)\leq \deg_{x_{c_{r+s}}}(x_ax_b) $. \\
	In particular, let $ [a] $ be a minor with $ a_1=1 $. From the definition of Hankel matrices, it is clear that, in the minor $ [a] $, there exists exactly one entry $ x_1 $. Therefore, for a given pair of minors $ [a] $ and $ [b] $, we have $ \deg_{x_1}(u)\leq \deg_{x_1}(x_ax_b)\leq 2 $ for all terms $ u $ in $ [a][b] $. Similar argument shows  $ \deg_{x_n}(u)\leq \deg_{x_n}(x_ax_b)\leq 2 $ for all terms $ u $ in $ [a][b] $.
\end{observe}

\begin{observe}\label{laplacian_expansion}
	Let $ [a] $ be a maximal minor of the  matrix $ \hmatTagSize{\sigma}{r} $. The Laplace expansion of $ [a] $ over the first row is
	$$[a]=\sum_{j=1}^{r}(-1)^{j+1}x_{a_j-j+1}[a_1+1,\ldots,a_{j-1}+1,\hat{a}_j,a_{j+1}\ldots,a_r].$$
	In particular 
	$$[a]=x_1[a_2,\ldots,a_r] + \tilde{H}$$
	where $\tilde{H}$ is the remaining factors of the Laplace expansion. Note for all terms $ u $ of $ \tilde{H} $, we have $ \deg_{x_1}(u)=0 $.\\
	The Laplace expansion of $ [a] $ over the last row is
	$$[a]=\sum_{j=1}^{r}(-1)^{j+1}x_{a_j+r-j}[a_1,\ldots,a_{j-1},\hat{a}_j,a_{j+1}-1\ldots,a_r-1].$$
	In particular 
	$$[a]=x_n[\chain{a}{r-1}] + H$$
	where $ H $ is the remaining factors of the Laplace expansion. Note for all terms $ u $ of $ H $, we have $ \deg_{x_n}(u)=0 $.

\end{observe}

\begin{exmp}
	Let $ n=10 $. Let $ [4,7,10] $ be a maximal minor in $ \hmatTagSize{(1,10)}{3} $. The Laplace expansion over the last row is
	$$ [4,7,10]=x_{10}[4,7]-x_8[4,9]+x_6[6,9]. $$ 
	Here, $ H=-x_8[4,9]+x_6[6,9] $ is the analogue of the one of \autoref{laplacian_expansion}. 
\end{exmp}

\begin{dfn}
	Let $ (\sigma,a) \geq_c (\gamma,b) $ be a pair of labeled chains. We say the product $ [a][b] $ has standard representation if 
	$$ [a][b]=\sum_i\lambda_i[c^{(i)}][d^{(i)}] $$
	such that  $ \lambda_i\in K $ and  $ (\sigma,c^{(i)}) \geq_c (\gamma,d^{(i)}) $ is a standard form of shape $ (\sigma,r) , (\gamma,s) $ for all $ i $. Moreover, 
	$$ \initial{[a][b]}>_{\tau} \initial{[c^{(i)}][d^{(i)}]} $$
	for all $ i>1 $. 
\end{dfn}

\begin{lem}\label{str_law}
	Let $ (\sigma,a) \geq_c (\gamma,b) $ be a pair of labeled chains. Then $ [a][b] $ has a standard representation.
	\begin{proof}
		One needs to repeat the following steps for finitely many times to obtain the standard representation of $ [a][b] $. One notes that this process eliminates the leading term in each repetition. Moreover, these terms are products of chains of lengths $ r $ and $ s $. Thus they are bounded from below with respect to $ \leq_\tau $. Hence, the algorithm converges. 
		\begin{enumerate}[Step(1)]
			\item Consider 
			$$ \delta= [a][b]-[c^{(1)}][d^{(1)}] $$ 
			where $ (\sigma,c^{(1)})\geq_c(\gamma,d^{(1)}) $ is the unique standard form obtained by applying \autoref{reduction_uniquness} on $ (\sigma,a)\geq_c(\gamma,b) $. From \autoref{anti_diagonal_reduction} one knows  $ (\sigma,c^{(1)})\geq_c(\gamma,d^{(1)}) $ has shape $ (\sigma,r),(\gamma,s) $. 
			In particular, we have $ \initial{\delta}<_\tau\initial{[a][b]} $. 
			\item Consider the pair of standard labeled chains $ (\sigma,c^{(2)})\geq_c(\gamma,d^{(2)}) $ by applying \autoref{well_def_remark} and \autoref{reduction_uniquness}. 
			\item Update $ \delta $ with $ \delta= \delta-\lambda_2[c^{(2)}][d^{(2)}] $ and return to step$ (2)$.
		\end{enumerate}
		
	\end{proof}
\end{lem}

\begin{lem}\label{well_def_remark}
	Let $ (\sigma,a) \geq_c (\gamma,b) $ be a pair of labeled chains and $ \delta=[a][b]-\sum_i\lambda_i[c^{(i)}][d^{(i)}] $ be obtained by finitely many times repeating the steps in \autoref{str_law}. Then $ \initial{\delta} $ admits a well-defined pair of labeled chains like $ (\sigma,e) \geq_c (\gamma,f) $ of shape $ (\sigma,r),(\gamma,s) $. 
	\begin{proof}
		We proceed by double induction on length of $ a=\chain{a}{r} $ and $ b=\chain{b}{s} $. When $ r=1 $ and $ s=1 $, it is trivial. Let $ r>1 $ and $ s>1 $ be positive integers. 
		The factors of $ \delta $ have the following  properties by construction:
		\begin{enumerate}
			\item  $ (\sigma,c^{(i)}) \geq_c (\gamma,d^{(i)}) $ is standard for all $ i $. Moreover, $ (\sigma,c^{(1)}) \geq_c (\gamma,d^{(1)}) $ is the standard form of $ (\sigma,a) \geq_c (\gamma,b) $. 
			\item  $ \initial{[c^{(j)}][d^{(j)}]}<_\tau\initial{[c^{(i)}][d^{(i)}]} $ for $ j>i $. In particular, $ \initial{[c^{(i)}][d^{(i)}]}<_\tau\initial{[a][b]} $ for all $ i>1 $. 
		\end{enumerate}
		By virtue of $ (1),(2) $ and \autoref{reduction_uniquness},  $ x_{c^{(i)}}x_{d^{(i)}} $'s are distinct for all $ i $. 
		From \autoref{G-basis_pruduct}, we have  chains $ e=\chain{e}{r} $ and $ f=\chain{f}{s} $ such that $ \lambda x_ex_f=\initial{\delta} $ for some $ \lambda\in K $. Note that from  
		\autoref{monomial_x_1_x_n_deg}, $ \deg_{x_{1}}(x_ax_b)\leq 2 $, $ \deg_{x_{1}}(x_{c^{(i)}}x_{d^{(i)}})\leq 2 $, $ \deg_{x_{n}}(x_ax_b)\leq 2 $ and $ \deg_{x_{n}}(x_{c^{(i)}}x_{d^{(i)}})\leq 2 $ for all $ i $.
		It is important to recall that $ \lambda x_ex_f $ is nevertheless some term of $ [a][b] $ or $ [c^{(i)}][d^{(i)}] $ for some $ i $. Thus \autoref{monomial_x_1_x_n_deg} implies $ \deg_{x_{1}}(x_ex_f)\leq 2 $ and $ \deg_{x_{n}}(x_ex_f)\leq 2 $.\\
		When  $ \deg_{x_{1}}(x_ax_b)=0 $ or $ \deg_{x_{1}}(x_ex_f)=0 $ or $ \sigma_1=\gamma_1 $, we always have $ \sigma_1\leq e_1 $ and $ \gamma_1\leq f_1 $. Note that, $ \deg_{x_{1}}(x_ax_b)=2 $, requires $ \sigma_1=\gamma_1 $ since $ (\sigma,a) \geq_c (\gamma,b) $ is a pair of labeled chains. So it falls into the previous case. Therefore,  to prove the well-definity on the left, it remains to consider the  case $ \sigma_1=1 $,  $ \gamma_1=2 $, $ \deg_{x_1}(x_ax_b)=1 $ (i.e $ a_1=1 $)  and $ \deg_{x_1}(x_ex_f)=1 $. Moreover, $ \deg_{x_{n}}(x_ax_b)=0 $ or $ \deg_{x_{n}}(x_ex_f)=0 $ or $ \sigma_2=\gamma_2 $ clearly implies $ e_r\leq \sigma_2 $ and $ f_s\leq \gamma_2 $. In particular, $ \deg_{x_{n}}(x_ax_b)=2 $ or $ \deg_{x_{n}}(x_ex_f)=2 $  requires $ \sigma_2=\gamma_2 $. Therefore it falls into the previous case. Hence it remains to consider the case  $ \deg_{x_{n}}(x_ax_b)= \deg_{x_{n}}(x_ex_f)=1 $  and $ \sigma_2\neq\gamma_2 $.

		To show the well-definity on the right, we split the rest of the proof with respect to value of $ \Delta $. 
		\begin{enumerate}[(I)]
			\item Suppose $ \Delta\{(\sigma,a) \geq_c (\gamma,b)\}=0 $. By definition of $ \Delta $, we have    $ a_r\leq \gamma_2 $ when $ r\leq s $ (or $ b_s\leq \sigma_2 $ when $ r> s $).   \\
			\begin{enumerate}[(i)]
				\item When $ r\geq s $, consider $ (e,f) $ as standard form in the sense of \autoref{ugly_remark}. From \autoref{ugly_remark} $ (II),(i) $, we have $ e_1=1 $. In particular, $ \deg_{x_{1}}(x_{e}x_{f})=1 $, yields  $ \sigma_1\leq e_1 $ and $ \gamma_1\leq f_1 $. Recall that $ \deg_{x_{n}}(x_{e}x_{f})\leq \deg_{x_{n}}(x_{a}x_{b}) $ from \autoref{monomial_x_1_x_n_deg}. Now, $ (e,f) $ being standard, the definition of $ \Delta $ and \autoref{ugly_remark} part $ (III) $ yields $ e_r\leq \sigma_2 $ and $ f_s\leq \gamma_2 $. 	
				\item When $ r < s $. By \autoref{laplacian_expansion}, we have 
				$$ \lambda x_ex_f=\initial{x_1[a_2,\ldots,a_{r}][b]-x_1\sum_i\lambda_i[c^{(i)}_2,\ldots,c^{(i)}_{r}][d^{(i)}]} $$
				where $ i $ runs through all $ c^{(i)} $'s with $ c^{(i)}_1=1 $.  Now, \autoref{G-basis_pruduct} admits the existence of  chains $ \tilde{e}=\tilde{e}_1,\ldots,\tilde{e}_{r-1} $ and $ \tilde{f}=\tilde{f}_1,\ldots,\tilde{f}_{s} $ such that 
				$$ \lambda x_{\tilde{e}}x_{\tilde{f}}=\lambda x_ex_f/x_1=\initial{[a_2,\ldots,a_{r}][b]-\sum_i\lambda_i[c^{(i)}_2,\ldots,c^{(i)}_{r}][d^{(i)}]}. $$
				Moreover, $ \deg_{x_2}(x_{\tilde{e}}x_{\tilde{f}})=\deg_{x_2}(x_{e}x_{f})\leq \deg_{x_2}(x_ax_b)\leq 1 $ by construction and \autoref{monomial_x_1_x_n_deg}. This implies that if $ \deg_{x_2}(x_{\tilde{e}}x_{\tilde{f}})=1 $, then $ 2 $ is the smallest entry in chains $ \tilde{e} $ and $ \tilde{f} $. Thus, from \autoref{ugly_remark} part $ (II),(i) $,  we can assume $ 2<\tilde{e}_1 $  by considering $ (\tilde{f},\tilde{e}) $ as standard form.  By reseting  notations, we have the chains $ e=1,\tilde{e}_1,\ldots,\tilde{e}_{r-1} $ and $ f=\tilde{f} $. It is clear that $ \sigma_1\leq e_1 $ and $ \gamma_1\leq f_1 $. 
				In particular, \autoref{ugly_remark} $ (III) $ and $ \deg_{x_n}(x_{\tilde{e}}x_{\tilde{f}})\leq \deg_{x_n}(x_{a}x_{b}) $ implies  $ e_{r}\leq \sigma_2 $ and $ f_s\leq \gamma_2 $. The reader notes that, we can not apply induction hypothesis on $ [a_2,\ldots,a_{r}][b]-\sum_i\lambda_i[c^{(i)}_2,\ldots,c^{(i)}_{r}][d^{(i)}] $, as $ (\sigma,c^{(i)}_2,\ldots,c^{(i)}_{r})\geq_c(\gamma,d^{(i)}) $'s are not necessarily standard in this context.  
			\end{enumerate}
			\item Suppose $ \Delta\{(\sigma,a) \geq_c (\gamma,b)\}=1 $. By definition of $ \Delta $, we have $ a_r=n $ and  $ \gamma_2=n-1 $ when $ r\leq s $  (or $ b_s=n $ and $ \sigma_2=n-1 $ when $ r> s $).\\
			\begin{enumerate}[(i)]
				\item When $ r\leq s $.  Form \autoref{laplacian_expansion}, we have 
				$$ \lambda x_ex_f=\initial{x_n[\chain{a}{r-1}][b]-x_n\sum_i[c^{(i)}_1\ldots,c^{(i)}_{r-1}][d^{(i)}]} $$
				where $ i $ runs through all $ c^{(i)} $'s with $ c^{(i)}_r=n $. By virtue of \autoref{G-basis_pruduct}, there exists chains $ \tilde{e}=\tilde{e}_1,\ldots,\tilde{e}_{r-1} $ and $ \tilde{f}=\tilde{f}_1,\ldots,\tilde{f}_{s} $ such that 
				\begin{equation}\label{eq_1}
				\lambda x_{\tilde{e}}x_{\tilde{f}}=\lambda x_ex_f/x_n=\initial{[\chain{a}{r-1}][b]-\sum_i[c^{(i)}_1\ldots,c^{(i)}_{r-1}][d^{(i)}]}.
				\end{equation}
				Note that  $ (\sigma,c^{(i)}_1,\ldots,c^{(i)}_{r-1})\geq_c(\gamma,d^{(i)}) $'s in \eqref{eq_1} are clearly standard forms here. Thus, the induction hypothesis implies  $ \sigma_1\leq \tilde{e}_1 $ and $ \gamma_1\leq \tilde{f}_1 $, in particular $ \tilde{e}_1=1 $. By applying \autoref{ugly_remark}, one shall consider $ (\tilde{e},\tilde{f}) $ as a standard form. One notes that \autoref{ugly_remark} part $ (I) $, implies that $ \tilde{e}_1 $ is not repositioned after considering $ (\tilde{e},\tilde{f}) $ as a standard form. Hence, one still has $ \tilde{e}_1=1 $.  On the other hand,  $ \deg_{x_{n-1}}(x_{\tilde{e}}x_{\tilde{f}})=\deg_{x_{n-1}}(x_{e}x_{f})\leq\deg_{x_{n-1}}(x_{a}x_{b})\leq 1 $ by construction and \autoref{monomial_x_1_x_n_deg}. From \autoref{ugly_remark} part $ (III) $, one has $ \tilde{e}_{r-1}<n-1 $ and $ \tilde{f}_{s}\leq n-1 $.  Hence, by reseting  notations, one has  chains $ e=\tilde{e}_1,\ldots,\tilde{e}_{r-1},n $ and $ f=\tilde{f}_1,\ldots,\tilde{f}_s $. Thus $ e_r\leq \sigma_2 $ and $ f_s\leq \gamma_2 $. In particular, $ \sigma_1\leq \tilde{e}_1 $ and $ \gamma_1\leq \tilde{f}_1 $ is a consequence of $ \tilde{e}=1 $. \\
				\item When $ r>s $. One can argue  similar to the last case. By definition of $ \Delta $, one has $ b_s=n $ and $ \sigma_2=n-1 $. From \autoref{laplacian_expansion}, one has 
				$$ \lambda x_ex_f=\initial{x_n[a][\chain{b}{s-1}]-x_n\sum_i[c^{(i)}][d^{(i)}_1\ldots,d^{(i)}_{s-1}]} $$
				where $ i $ runs through all $ d^{(i)} $'s with $ d^{(i)}_s=n $. By virtue of \autoref{G-basis_pruduct}, there exists chains $ \tilde{e}=\tilde{e}_1,\ldots,\tilde{e}_r $ and $ \tilde{f}=\tilde{f}_1,\ldots,\tilde{f}_{s-1} $ such that 
				\begin{equation}\label{eq_2}
				\lambda x_{\tilde{e}}x_{\tilde{f}}=\lambda x_ex_f/x_n=\initial{[a][\chain{b}{s-1}]-\sum_i[c^{(i)}][d^{(i)}_1\ldots,d^{(i)}_{s-1}]}.
				\end{equation}
				Note that  $ (\sigma,c^{(i)})\geq_c(\gamma,d^{(i)}_1,\ldots,d^{(i)}_{s-1}) $'s in \eqref{eq_2} are clearly standard forms here. Thus, the induction hypothesis implies  $ \sigma_1\leq \tilde{e}_1 $ and $ \gamma_1\leq \tilde{f}_1 $, in particular $ \tilde{e}=1 $.
				By applying \autoref{ugly_remark}, one shall consider $ (\tilde{e},\tilde{f}) $ as a standard form. Note that \autoref{ugly_remark} part $ (II),(i) $ implies $ \tilde{e}_1=1
				$.  Hence, one still has $ \tilde{e}_1=1 $ in particular  $ \sigma_1\leq \tilde{e}_1 $ and $ \gamma_1\leq \tilde{f}_1 $.  Recall that $ \deg_{x_{n-1}}(x_{\tilde{e}}x_{\tilde{f}})=\deg_{x_{n-1}}(x_{e}x_{f})\leq\deg_{x_{n-1}}(x_{a}x_{b})\leq 1 $ by construction and \autoref{monomial_x_1_x_n_deg}. From \autoref{ugly_remark} part $ (III) $, one has $ \tilde{e}_{r-1}\leq n-1 $ and $ \tilde{f}_{s}< n-1 $.  Hence, by reseting  notations, one has  chains $ e=\tilde{e}_1,\ldots,\tilde{e}_{r} $ and $ f=\tilde{f}_1,\ldots,\tilde{f}_{s-1},n $. Thus $ e_r\leq \sigma_2 $ and $ f_s\leq \gamma_2 $.\\
			\end{enumerate}
		\end{enumerate}
	\end{proof}
\end{lem}

\begin{theorem}\label{main_theorem}
	The family $ \f=\f^{(1,n)}\cup\f^{(1,n-1)}\cup\f^{(2,n)}\cup\f^{(2,n-1)} $ has the following features:
	\begin{enumerate}[(1)]
		\item Every product $ \prod_{(\sigma,a)}\detIdeal{\sigma}{a} $ of ideals in $ \f $ has linear resolution.
		\item Computing the initial ideals commutes over products  $ \initial{\prod_{(\sigma,a)}\detIdeal{\sigma}{a}}=\prod_{(\sigma,a)}\initial{\detIdeal{\sigma}{a}} $, in particular the natural generators form a Gr\"{o}bner basis. 
		\item The multi-Rees algebra $ \rees{\detIdeal{\sigma}{a}:\detIdeal{\sigma}{a}\in\f} $ is defined by a quadratic Gr\"{o}bner basis with respect to $ \leq_{\alpha} $, it is Koszul,  normal, Cohen-Macaulay domain. Moreover, the natural algebra generators form a Sagbi basis.
	\end{enumerate}  
	\begin{proof}
		 It is enough to prove $ (3) $. We apply \cite[Crollary 2.2]{Conca_Herzog_Valla}. The binomials of the form $ x_i z_{\sigma,a}-x_j z_{\sigma,c}\ $ and $ z_{\sigma,a} z_{\sigma,b}-z_{\sigma,c} z_{\sigma,d} $ in $ G $  lifts to $ \kerRees $ by \cite{NAM_PAPER}. The binomials $ z_{\sigma,a} z_{\gamma,b}-z_{\sigma,c} z_{\gamma,d} $ in $ G $ lifts to $ z_{\sigma,a} z_{\gamma,b}-\sum_iz_{\sigma,c^{(i)}} z_{\gamma,d^{(i)}} $ in $ \kerRees $ where the indices are obtained from  \autoref{str_law}. This admits a quadratic Gr\"{o}bner basis for   $ \kerRees $. Therefor,  $ \reesSymbol $ is Koszul. In particular, by virtue of \cite[Preposition 1.1]{Conca_Herzog_Valla}, the algebraic generators of $ \reesSymbol $ form a Sagbi basis, which is equivalent to $ (2) $. Now, $ (1) $ is  consequence of Koszulness of $ \reesSymbol $ and the multi-graded version of the theorem of Blum \cite{Blum}. \\
		 To prove that $ \reesSymbol $ is normal, Cohen-Macaulay domain, by \cite[Corollary 2.3]{Conca_Herzog_Valla}, is is enough to prove that $ \reesSymbol^{\operatorname{in}} $ is normal. Recall that the term order $ \leq_{\alpha} $ picks non-standard monomials as the leading terms of the elements in $ G $. Moreover, every non-square free monomial of degree two in indeterminates $ \bar{z} $ is standard. Thus $ \operatorname{in}_{\leq_{\alpha}}(\kerReesIn) $ is  square free. Hence \cite[Prposition 13.15]{Sturmfels_marked_polynomials} yields the normality of $ \reesSymbol^{\operatorname{in}} $.
	\end{proof}	
\end{theorem}

We conclude this paper by explaining why the family of close cuts of Hankel matrices are interesting. 

\begin{rem}\label{importance}
 	Let $ x_i,\ldots,x_j $ be an interval of indeterminates of $ S $ where $ i\leq j $. Let  $ \hmatTagSize{(i,j)}{} $ and $ \f^{(i,j)} $ be defined similar to $ \hmatTagSize{(1,n)}{} $ and $ \f^{(1,n)} $. Let $ \f=\cup_{i\leq j}\f^{(i,j)} $.
 	 We expect \autoref{main_theorem}, $ (1) $ to extend for $ \f $. As we saw, one standard approach is via Sagbi deformations.  However, it is easy to see that this is not the case for  \autoref{main_theorem}, $ (2) $. For $ n\geq 6 $, we have  $ \initial{\detIdeal{(1,n)}{2}\detIdeal{(3,n)}{2}}\neq\initial{\detIdeal{(1,n)}{2}}\initial{\detIdeal{(3,n)}{2}} $ which is equivalent to $ \reesInitial{\initial{I}:I\in\f}\neq\initial{\rees{I:I\in\f}} $. Hence, the kernel of $ \reesInitial{\initial{I}:I\in\f} $ does not lift to  the kernel of $ \rees{I:I\in\f} $ (see \cite[Proposition 1.1]{Conca_Herzog_Valla}). Therefore,  Sagbi deformation method fails. Nevertheless, We still expect \autoref{main_theorem_initial} to extend for $ \f $.\\ 
 	 Conca and Nam, in separated papers, prove that the product $ I=I_1\ldots I_l $ of ideals in $ \f^{(1,n)} $  have a nice primary decomposition   given by intersection of symbolic powers of ideals in $ \f^{(1,n)} $ containing $ I $ (see \cite[Theorem 3.12]{CONCA_STR_LAW} and \cite[Theorem 3.25]{NAM_PAPER}). The author refers to standard text books in commutative algebra for the definition of symbolic powers. Similar feature is provided for generic matrices, however, with ordinary powers by  Berget, Bruns and Conca (see \cite[Corollary 2.3]{Berget_Bruns_Conca} and \cite[Theorem 3.4]{Bruns_Conca}). This might raise the question that whether this behavior is expected for $  \f $. Unfortunately, this is not the case.   Consider $ I=\detIdeal{(1,6)}{2}\detIdeal{(3,6)}{2} $. We have $ \operatorname{Ass}(I)=\{\detIdeal{(1,6)}{2},\detIdeal{(3,6)}{2},\detIdeal{(2,6)}{1}\} $, where $ \operatorname{Ass}(I) $ is the associated primes of $ I $. Thanks to  \cite[Theorem 3.8]{CONCA_STR_LAW}, one can check  $ I\subset\detIdeal{(3,6)}{2} $, $ I\subset\detIdeal{(1,6)}{2}^{(2)} $ and $ I\subset\detIdeal{(2,6)}{1}^{(3)} $. The inclusions is strict and increasing the symbolic exponent will defy the inclusion. In particular,  $ I\subsetneq\detIdeal{(3,6)}{2}\cap \detIdeal{(1,6)}{2}^{(2)}\cap\detIdeal{(2,6)}{1}^{(3)} $. Hence, we can not expect a similar behavior of the primary decompositions for $ \f $. Nevertheless, a nice primary decomposition is expected for $ \f^{(1,n)}\cup\f^{(1,n-1)}\cup\f^{(2,n)}\cup\f^{(2,n-1)} $.

 \end{rem}

\bibliographystyle{elsarticle-num} 
\bibliography{references}
\end{document}